\documentclass[a4paper]{article}
\usepackage [a4paper,left=3 cm,bottom=2.5cm,right=3cm,top=2.5cm]{geometry}
\usepackage[english]{babel}
\usepackage{amssymb}
\usepackage{amsmath,amsthm}
\usepackage{subcaption}
\usepackage{tabularx,array}
\usepackage{graphicx}
\usepackage{enumitem} 

\newtheorem{theorem}{Theorem}
\newtheorem{lemma}{Lemma}
\newtheorem{corollary}{Corollary}
\newtheorem{proposition}{Proposition}
\newtheorem{definition}{Definition}

\newcommand{\keywords}[1]{\par\addvspace\baselineskip
\noindent\enspace\ignorespaces#1}
\begin{document}

\title{Invariant density adaptive estimation for ergodic jump diffusion processes over anisotropic classes.}
\author{Chiara Amorino$^{*}$, Arnaud Gloter\thanks{Laboratoire de Math\'ematiques et Mod\'elisation d'Evry, CNRS, Univ Evry, Universit\'e Paris-Saclay, 91037, Evry, France.}} 

\maketitle

\begin{abstract}
We consider the solution $X=(X_t)_{t \ge 0}$ of a multivariate stochastic
differential equation with Levy-type jumps and with unique invariant probability measure with density $\mu$. We assume that a continuous record of observations $X^T = (X_t)_{0 \le t \le T}$ is available. \\ 
In the case without jumps, Reiss and Dalalyan \cite{RD} and Strauch \cite{Strauch} have found convergence rates of invariant density estimators, under respectively isotropic and anisotropic H{\"o}lder smoothness constraints, which are considerably faster than those known from standard multivariate density estimation. \\
We extend the previous works by obtaining, in presence of jumps, some estimators which have the same convergence rates they had in the case without jumps for $d \ge 2$ and a rate which depends on the degree of the jumps in the one-dimensional setting. \\
We propose moreover a data driven bandwidth selection procedure based on the Goldenshluger and Lepski method \cite{Adaptive} which leads us to an adaptive non-parametric kernel estimator of the stationary density $\mu$ of the jump diffusion $X$.

\end{abstract}

\section{Introduction}
Diffusion phenomena arise from a Markovian stochastic modeling and as a solution of SDEs with or without jumps in many areas of applied mathematics. Their investigation concerns different mathematical branches and therefore research interest in questions such as existence and regularity of solutions of stochastic differential equations has constantly grown over the past years. \\
The study of the statistical properties of diffusion models has emerged since such models are widely used for applications in finance and biology. Diffusion processes with jumps, in particular, have been used in neuroscience for instance in \cite{Neuro} while in finance they have been introduced to model the dynamic of asset prices \cite{Kou}, \cite{Merton}, exchange rates \cite{Bates}, or volatility processes \cite{BarShe}. 

In this work, we aim at estimating adaptively the invariant density $\mu$ associated to the process $(X_t)_{t \ge 0}$, solution of the following multivariate stochastic
differential equation with Levy-type jumps:
\begin{equation}
X_t= X_0 + \int_0^t b( X_s)ds + \int_0^t a(X_s)dW_s + \int_0^t \int_{\mathbb{R}^d \backslash \left \{0 \right \} }
\gamma(X_{s^-})z \tilde{\mu}(ds,dz),
\label{eq: model intro}
\end{equation}
where $W$ is a $d$-dimensional Brownian motion and $\bar{\mu}$ a compensated Poisson random measure with a possible infinite jump activity. We assume that a continuous record of observations $X^T = (X_t)_{0 \le t \le T}$ is available. 

Practical concerns raise new questions such as the dependence of statistical features on the observation scheme: it is, for the applications, a subject of interest
to consider basic questions in different observation scenarios.
From a theoretical point of view, it is however also of substantial interest to work under the assumption that a continuous record of the diffusion considered is available. \\
In this framework, it belongs to the folklore of the statistics for stochastic processes without jumps that the invariant density can be estimated under standard nonparametric assumptions with a parametric rate (cfr Chapter 4.2 in \cite{Kut}). The proof relies on the existence of diffusion local time and its properties and so such a result is restricted to the one - dimensional setting. 

Regarding the literature on statistical properties of multidimensional diffusion processes in the continuous case, an important reference is given by Reiss and Dalalyan in \cite{RD}, where they show an asymptotic statistical equivalence for inference on the drift in the multidimensional diffusion case. As a by-product of the study they prove, under isotropic H{\"o}lder smoothness constraints, convergence rates of invariant density estimators for pointwise estimation which are faster than those known from standard multivariate density estimation. Their result relies on upper bounds on the variance of additive diffusion functionals, proven by an application of the spectral gap inequality in combination with a bound on the transition density of the process. \\
Still in the continuous case, in a recent paper, Strauch \cite{Strauch} has extended their work by building adaptive estimators in the multidimensional diffusion case which achieve fast rates of convergence over anisotropic H{\"o}lder balls. \\
The notion of anisotropy plays an important role. Indeed, the smoothness properties of elements of a function space may depend on the chosen direction of $\mathbb{R}^d$. \\
The Russian school considered anisotropic spaces from the beginning of the theory of function spaces in 1950-1960s (in \cite{Nik} the author takes account of the developments). However, results on minimax rates of convergence in classical statistical models were rare for a lot of time. \\
The question of optimal bandwidth selection based on i. i. d. observations for density estimation with respect to sup - norm risk was not completely solved until the pretty recent developments gathered in \cite{Lep}. The methodology detailed in Goldenshluger and Lepski \cite{Adaptive} inspired the data-driven selection
procedure of the bandwidth of the kernel estimator proposed by many authors such as Strauch in \cite{Strauch} and Comte, Prieur and Samson in \cite{Main adapt} and provides the starting point for the study of our adaptive procedure as well.

In this paper, we provide a non-parametric estimator of the invariant density $\mu$ with a fully data-driven
procedure of the bandwidth. 
We  propose  to  estimate  the  invariant  density $\mu$ by  means  of  a  kernel  estimator, we therefore introduce some kernel function $K: \mathbb{R} \rightarrow \mathbb{R}$. A natural estimator of $\mu$ at $x \in \mathbb{R}^d$ in the anisotropic context is given by
$$\hat{\mu}_{h,T}(x) = \frac{1}{T \prod_{l = 1}^d h_l} \int_0^T \prod_{m = 1}^d K(\frac{x_m - X_u^m}{h_m}) du,$$
where $h = (h_1, ... , h_d)$ is a multi - index bandwidth, which will be chosen through the data-driven selection
procedure. 
We first prove some bounds on the transition semigroup and on the transition density that will be useful to find sharp upper bounds on the variance of integral functionals of the diffusion $X$. Through them, we find the following convergence rates for the pointwise estimation of the invariant density of our diffusion with jumps:
$$\mathbb{E} [|\hat{\mu}_{h,T}(x) - \mu (x)|^2] \underset{\sim}{<}
\begin{cases}
\frac{(\log T)^{(2 - \frac{(1 + \alpha)}{2}) \lor 1}}{T} \qquad \mbox{for } d = 1, \\
\frac{\log T}{T} \qquad \mbox{for } d = 2, \\
T^{- \frac{2\bar{\beta}}{2\bar{\beta}+ d - 2}} \qquad \mbox{for } d \ge 3,
\end{cases}
$$
where $\alpha \in (0, 2)$ is the degree of jumps activity of the L\'evy process and $\bar{\beta}$ is the harmonic mean smoothness of the invariant density over the $d$ different dimensions. \\
We remark that the rate we find for $d \ge 3$ is the same Strauch found in \cite{Strauch} in absence of jumps, which is also the rate gathered in \cite{RD} up to replacing the mean smoothness with $\beta$, the common smoothness over the $d$ dimensions. \\
The case $d = 1$ evidences the main difference between what happens with and without jumps. Indeed, if in the continuous case the optimal convergence rate was $\frac{1}{T}$, now the rate we found is between $\frac{\log T}{T}$ and $\frac{(\log T)^\frac{3}{2}}{T}$. It is worth noting here that such a convergence rate is not necessarily the optimal one in the jumps framework. As a matter of fact in the continuous case different approaches, as the diffusion local time, have been used to get the rate $\frac{1}{T}$; we do not exclude the possibility that also in presence of jumps the implementation of other methods could lead to a convergence rate faster than the one presented here above for the mono-dimensional setting. \\
To complete the comparison to the continuous framework, we recall that in both \cite{RD} and \cite{Strauch} the convergence rate found in the case $d = 2$ was $\frac{(\log T)^4}{T}$ and so the convergence of the estimator seems being faster in presence of jumps than without them. The reason why it happens is that, to find the convergence rate, the transition density $(p_t)_{t \in \mathbb{R}^+}$ is needed to be upper bounded. If in \cite{RD} the authors assume to have $p_t(x,y) \le c (t^{- \frac{d}{2}} + t^{ \frac{3d}{2}})$ and in \cite{Strauch} Nash and Poincar\'e inequalities lead Strauch to a bound analogous to the one presented in \cite{RD}; Lemma \ref{lemma: bound transition} below provides us a different bound which guides us to a different rate. However, in absence of the term $t^{ \frac{3d}{2}}$ in the assumption before, which is the case for example considering a bounded drift, also in the continuous setting the convergence rate turns out being, as in the jump -diffusion case, equal to $\frac{\log T}{T}$. \\ 
It is moreover worth noting here that, if in \cite{RD} and \cite{Strauch} they needed to assume the existence of the transition density and a bound on it, we derive them through Lemma \ref{lemma: bound transition}: all the assumptions we need are directly on the model \eqref{eq: model intro}. \\
We no longer need to assume that the drift is of the form $b = - \nabla V$ (where $V \in \mathcal{C}^2$ is referred to as potential) as it was in both \cite{RD} and \cite{Strauch}.

After having provided the rates of
convergence of the estimators we finally propose, in the case $d \ge 3$, a fully data-driven selection
procedure of the bandwidth of the kernel estimator, inspired by the methodology
detailed in Goldenshluger and Lespki \cite{Adaptive}. The method has the decisive advantage of being
anisotropic: the bandwidths selected in each direction are in general different, which is coherent
with the possibly different regularities with respect to each variable. Finally, we prove that for the selected optimal bandwidth the following estimation holds:
\begin{equation}
\mathbb{E}[\left \| \hat{\mu}_{\tilde{h}} -  \mu \right \|^2_A] \le c_1 \inf_{h \in \mathcal{H}_T} (B(h) + V(h)) + c_1 e^{ - c_2 (\log T)^2},
\label{eq: tesi adatt}
\end{equation}
where we have denoted as $\left \| \cdot \right \|_A$ the $L^2$ norm on $A$, a compact subset of $\mathbb{R}^d$ and as $\mathcal{H}_T$ the set of candidate bandwidths; $B(h)$ is a bias term and $V(h)$ an estimate of the variance bound. We remark that the estimator leads to an automatic trade - off between the bias and the variance: the second term on the right hand side of \eqref{eq: tesi adatt} is indeed negligible compared to the first one. \\
Moreover, as the rate optimal choice $h(T)$ belongs to the set of candidate bandwidths $\mathcal{H}_T$, \eqref{eq: tesi adatt} turns out being
$$\mathbb{E}[\left \| \hat{\mu}_{\tilde{h}} -  \mu \right \|^2_A] \le c_1 T^{- \frac{2\bar{\beta}}{2\bar{\beta}+ d - 2}} + c_1 e^{ - c_2 (\log T)^2},$$
where $\bar{\beta}$ is the mean smoothness of the invariant density.

The paper is organised as follows. We give in Section 2 the assumptions on the process $X$. In section \ref{S: Construction_estimator} we define the anisotropic H{\"o}lder balls and we construct our estimator. Section \ref{S: Main_results} is devoted to the statements of our main results; which will be proven in the two following sections. In particular, we show how we get the convergence rates for the invariant density estimation in Section \ref{S: Proof_Main} while in Section \ref{S: proof adaptive} we prove the estimator we find through our bandwidth selection procedure is adaptive. Some technical results are moreover presented in the Appendix.

\section{Model Assumptions}
We consider the question of nonparametric estimation of the invariant density of a d-dimensional diffusion process X, assuming that a continuous record $X^T = \left \{  X_t, 0 \le t \le T \right \}$ up to time T is observed. This diffusion is given as a strong solution of the following stochastic differential equations with jumps:
\begin{equation}
X_t= X_0 + \int_0^t b( X_s)ds + \int_0^t a(X_s)dW_s + \int_0^t \int_{\mathbb{R}^d \backslash \left \{0 \right \} }
\gamma(X_{s^-})z \tilde{\mu}(ds,dz), \quad t \in [0,T], 
\label{eq: model}
\end{equation}
where $b : \mathbb{R}^d \rightarrow \mathbb{R}^d$, $a : \mathbb{R}^d \rightarrow \mathbb{R}^d \times \mathbb{R}^d$ and $\gamma : \mathbb{R}^d \rightarrow \mathbb{R}^d \times \mathbb{R}^d$; $W = (W_t, t \ge 0)$ is a d-dimensional Brownian motion and $\mu$ is a Poisson random measure on $[0, T] \times \mathbb{R}^d$ associated to the L\'evy process $L=(L_t)_{t \in [0,T]}$, with $L_t:= \int_0^t \int_{\mathbb{R}^d} z \tilde{\mu} (ds, dz)$. The compensated measure is $\tilde{\mu}= \mu - \bar{\mu}$; we suppose that the compensator has the following form: $\bar{\mu}(dt,dz): = F(dz) dt $, where conditions on the Levy measure $F$ will be given later. \\ 
The initial condition $X_0$, $W$ and $L$ are independent. \\ \\
In what follows, we suppose the following assumptions hold: \\ \\
\textbf{A1}: \textit{The functions $b(x)$, $\gamma(x)$ and $a(x)$ are globally Lipschitz and, for some $c \ge 1$,
$$c^{-1} \mathbb{I}_{d \times d} \le a(x) \le c \mathbb{I}_{d \times d}, $$
where $\mathbb{I}_{d \times d}$ denotes the $d \times d$ identity matrix. \\
Denoting with $|.|$ and $<., . >$ respectively the Euclidian norm and the scalar product in $\mathbb{R}^d$, we suppose moreover that there exists a constant $c > 0$ such that, $\forall x \in \mathbb{R}^d$, $|b(x)| \le c$.} \\
\\
Under Assumption 1 the equation (\ref{eq: model}) admits a unique non-explosive c\`adl\`ag adapted solution possessing the strong Markov property, cf \cite{Applebaum} (Theorems 6.2.9. and 6.4.6.). \\
\\
\textbf{A2 (Drift condition) }: \textit{ \\
 There exist $\tilde{C} > 0$ and $\tilde{\rho} > 0$ such that $<x, b(x)>\, \le -\tilde{C}|x|$, $\forall x : |x| \ge \tilde{\rho}$.
 } \\
\\
We furthermore need the following assumptions on the jumps: \\
\\
\textbf{A3 (Jumps) }: \textit{1.The L\'evy measure $F$ is absolutely continuous with respect to the Lebesgue measure and we denote $F(z) = \frac{F(dz)}{dz}$. \\
2. We suppose that there exist $c > 0$ such that for all $z \in \mathbb{R}^d$, $F(z) \le \frac{c}{|z|^{d + \alpha}}$, with $\alpha \in (0,2)$ and that $supp(F) = \mathbb{R}^d$. \\
3. The jump coefficient $\gamma$ is upper bounded, i.e. $\sup_{x \in \mathbb{R}^d}|\gamma(x)| := \gamma_{max} < \infty$. We suppose moreover that, $\forall x \in \mathbb{R}^d$, $Det(\gamma(x)) \neq 0$.\\
4. If $\alpha =1$, we require for any $0 < r < R < \infty$ $\int_{r <| z |< R} z F(z) dz =0$. \\
5. There exists $\epsilon > 0$ and a constant $c > 0$ such that $\int_{\mathbb{R}^d}|z|^2 e^{\epsilon |z|} F(z) dz \le c $.} \\
\\
As we will see in Lemma \ref{lemma: beta mixing} below, Assumption 2 ensures, together with the last points of Assumption 3, the existence of a Lyapunov function. The process $X$ admits therefore a unique invariant distribution $\pi$ and the ergodic theorem holds. \\
We assume the invariant probability measure $\pi$ of $X$ being absolutely continuous with respect to the Lebesgue measure and from now on we will denote its density as $\mu$: $ d\pi = \mu dx$. \\
For any set $S \subset \mathbb{R}^d$ we define $\mu (S) := \int_S \mu(x) dx$ and, by abuse of notation, we will write $\mu(f) : = \mathbb{E}[f (X_0)] = \int_{\mathbb{R}^d} f(x) \mu (x) dx$ for functions $f : \mathbb{R}^d \rightarrow \mathbb{R}$. \\
We define moreover $L^2 (\mu) := \left \{ f : \mathbb{R}^d \rightarrow \mathbb{R} : \int_{\mathbb{R}^d} |f(x)|^2 \mu (x) dx < \infty \right \}$ and \\
$L^1 (\mu) := \left \{ f : \mathbb{R}^d \rightarrow \mathbb{R} : \int_{\mathbb{R}^d} |f(x)| \mu (x) dx < \infty \right \}$. \\
For each $g\in L^1 (\mu)$ we denote as $\left \| g \right \|_{L^1 (\mu)} := \mu (|g|) $ the $L^1$ norm with respect to $\mu$ on $\mathbb{R}^d$. \\
The transition semigroup of the process $X$ on $L^1 (\mu)$ is $P_{t} f(x) := \mathbb{E} [f(X_t) | X_0 = x ]$. \\
The transition density is denoted by $p_{t}$ and it is such that $P_{t} f(x) = \int_{\mathbb{R}^d} f(y) p_{t} (x,y) dy$; we will see in Lemma \ref{lemma: bound transition} that it exists. \\ \\
The process $X$ is called $\beta$ - mixing if $\beta_X (t) = o(1)$ for $t \rightarrow \infty$ and exponentially $\beta$ - mixing if there exists a constant $\gamma > 0$ such that $\beta_X (t) = O(e^{- \gamma t})$ for $t \rightarrow \infty$, where $\beta_X$ is the $\beta$ - mixing coefficient of the process $X$ as defined in Section 1.3.2 of \cite{Mixing}. We recall that, for a Markov process $X$ with transition semigroup $(P_t)_{t \in \mathbb{R}^+}$ and $\mathcal{L}(X_0) = \eta$, the $\beta$ - mixing coefficient of $X$ is given by 
\begin{equation}
\beta_X (t) := \sup_{s \in \mathbb{R}^+} \int_{\mathbb{R}^d} \left \| P_t (x, .) - \eta P_{s + t} (x, .) \right \| \eta P_s (dx, .),
\label{eq: def beta}    
\end{equation}
where $\eta P_t = \mathcal{L}(X_t)$ and $\left \| \lambda \right \|$ stands for the total variation norm of a signed measure $\lambda$. \\
For the exponential mixing property of general multidimensional diffusions, the reader may consult Theorem 3 of Kusuoka and Yoshida \cite{Kus_Yos} for the $\alpha$ - mixing; Meyn and Tweedie \cite{May_Twe}, Stramer and Tweedie \cite{Str_Twe} and Veretennikov \cite{Ver} for the $\beta$ - mixing. The mixing property for general diffusions with jumps has been investigated by Masuda in \cite{18 GLM}.  \\
Now we mention the notion of exponential ergodicity in the sense of \cite{May_Twe}. \\
\begin{definition}
We say that $X$ is exponentially ergodic if it admits a unique invariant distribution $\pi$ and additionally if there exist positive constants c and $\rho$ for which, for each $f$ centered under $\mu$, 
$$\left \| P_{t} f \right \|_{L^1 (\mu)} \le c e^{- \rho t} \left \| f \right \|_{\infty}. $$
\label{def: exp ergodic}
\end{definition}
We will see in Lemma \ref{lemma: beta mixing} that both the exponential ergodicity and the exponential $\beta$ - mixing can be derived from our assumptions. \\
\\
In Lemmas \ref{lemma: beta mixing} and \ref{lemma: bound transition} below we will prove some bounds on the transition semigroup and on the transition density that will be useful to establish tight upper bounds on the variance
$$Var (\int_0^T f(X_s) ds), \qquad f \in L^2(\mu)$$
of integral functionals of the diffusion $X$.  \\
Bounds of this type were proven before, in \cite{RD} (cf. their Proposition 1), by combining estimates based on the spectral gap inequality and on upper bounds on the transition densities of $X$. 
Through them they prove, under isotropic H{\"o}lder smoothness constraints, convergence rates of invariant density estimators for pointwise estimation  which are considerably faster than those known from standard multivariate density estimation. \\
We replace the spectral gap inequality with a control from $L^1$ to $L^\infty$ given by the exponential ergodicity. Moreover, contrary to \cite{RD}, we don't need to assume that such controls hold true since we get them as consequence of Lemma \ref{lemma: bound transition} and \ref{lemma: beta mixing} below, having required some assumptions only directly on the model \eqref{eq: model}.  \\
In the next section we will construct adaptive estimators for the density in the multidimensional diffusion case with jumps, which achieve ‘fast’ rates of convergence over anisotropic H{\"o}lder balls.

\section{Construction of the estimator}\label{S: Construction_estimator}
In several cases, the regularity of some function $g: \mathbb{R}^d \rightarrow \mathbb{R}$ depends on the direction in $\mathbb{R}^d$ chosen. We thus work under the following anisotropic smoothness constraints.
\begin{definition}
Let $\beta = (\beta_1, ... , \beta_d)$, $\beta_i > 0$, $\mathcal{L} =(\mathcal{L}_1, ... , \mathcal{L}_d)$, $\mathcal{L}_i > 0$. A function $g : \mathbb{R}^d \rightarrow \mathbb{R}$ is said to belong to the anisotropic H{\"o}lder class $\mathcal{H}_d (\beta, \mathcal{L})$ of functions if, for all $i \in \left \{ 1, ... , d \right \}$,
$$\left \| D_i^k g \right \|_\infty \le \mathcal{L}_i \qquad \forall k = 0,1, ... , \lfloor \beta_i \rfloor, $$
$$\left \| D_i^{\lfloor \beta_i \rfloor} g(. + t e_i) - D_i^{\lfloor \beta_i \rfloor} g(.) \right \|_\infty \le \mathcal{L}_i |t|^{\beta_i - \lfloor \beta_i \rfloor} \qquad \forall t \in \mathbb{R},$$
for $D_i^k g$ denoting the $k$-th order partial derivative of $g$ with respect to the $i$-th component, $\lfloor \beta_i \rfloor$ denoting the largest integer strictly smaller than $\beta_i$ and $e_1, ... , e_d$ denoting the canonical basis in $\mathbb{R}^d$.
\end{definition}
From now on we deal with the estimation of the density $\mu$ belonging to the anisotropic H{\"o}lder class $\mathcal{H}_d (\beta, \mathcal{L})$. \\
\\
Given the observation $X^T$ of a diffusion $X$, solution of \eqref{eq: model}, we propose to estimate the invariant density $\mu$ by means of a kernel estimator.
To estimate some $\mu \in \mathcal{H}_d (\beta, \mathcal{L})$ we therefore introduce some kernel function $K: \mathbb{R} \rightarrow \mathbb{R}$ satisfying 
$$\int_\mathbb{R} K(x) dx = 1, \quad \left \| K \right \|_\infty < \infty, \quad \mbox{supp}(K) \subset [-1, 1], \quad \int_\mathbb{R} K(x) x^l dx = 0,$$
for all $l \in \left \{ 0, ... , M \right \}$ with $M \ge \max_i \beta_i$. \\
Denoting by $X_t^j$, $j \in \left \{ 1, ... , d \right \}$ the $j$-th component of $X_t$, $t \ge 0$, a natural estimator of $\mu$ at $x= (x_1, ... , x_d)^T \in \mathbb{R}^d$ in the anisotropic context is given by 
\begin{equation}
\hat{\mu}_{h,T}(x) = \frac{1}{T \prod_{l = 1}^d h_l} \int_0^T \prod_{m = 1}^d K(\frac{x_m - X_u^m}{h_m}) du.
\label{eq: def estimator}
\end{equation}
As we will see in Section \ref{S: Adaptive}, a main question concerns the choice of the multi-index bandwidth $h = (h_1, ... , h_d)^T$.

\section{Main results}\label{S: Main_results}
\subsection{Convergence rates for invariant density estimation}
We want to investigate on the convergence rates for invariant density estimation. In order to determine the asymptotic behaviour of our estimator for $T \rightarrow \infty$, we study the variance of general additive functionals of $X$ in $d$ dimension. To do so, we need some properties as the exponential ergodicity of the process and a bound on the transition density. Such properties will be derived from our assumptions through the following lemmas, that we will prove in the appendix. \\
The following bounds on the transition density and on the transition semigroup hold true.
\begin{lemma}
Suppose that A1 - A3 hold. Then, for $T \ge 0$, there exists a transition density $p_{t} (x,y)$ for which for any $t \in [0, T]$ there are a $c_0 > 0$ and a $\lambda_0 > 0$ such that, for any pair of points $x, y \in \mathbb{R}^d$, we have
$$p_{t} (x,y) \le c_0 (t^{- \frac{d}{2}} e^{- \lambda_0 \frac{|y - x|^2}{t}} + \frac{t}{(t^\frac{1}{2} + |y - x|)^{d + \alpha}}).$$
\label{lemma: bound transition}
\end{lemma}
\begin{lemma}
Suppose that A1 - A3 hold. Then the process $X$ is exponentially ergodic and exponentially $\beta$ - mixing.
\label{lemma: beta mixing}
\end{lemma}
On the basis of the two previous lemmas we can prove the following bound on the variance, which is the heart of the study on the convergence rate.
\begin{proposition}
Suppose that A1 - A3 hold and let $f: \mathbb{R}^d \rightarrow \mathbb{R}$ be a bounded, measurable function with support $\mathcal{S}$ satisfying $|\mathcal{S}| < 1$. Then, there exists a constant C independent of $f$ such that
\begin{itemize}
    \item[$\bullet$] $Var (\int_0^T f (X_t) dt ) \le C T \left \| f \right \|_\infty^2 |\mathcal{S}|^2(1 + (\log(\frac{1}{|\mathcal{S}|}))^{2 - \frac{(1 + \alpha)}{2}} + \log(\frac{1}{|\mathcal{S}|}))\quad$ for $d=1$,
    \item[$\bullet$] $Var (\int_0^T f (X_t) dt ) \le C T\left \| f \right \|^2_{\infty}|\mathcal{S}|^2 (1 + \log(\frac{1}{|\mathcal{S}|})) \quad$ for $d=2$,
    \item[$\bullet$] $Var (\int_0^T f (X_t) dt ) \le C T \left \| f \right \|_\infty^2 |\mathcal{S}|^{1 + \frac{2}{d}} \quad$ for $d \ge 3$.
\end{itemize}
\label{prop: variance bound}
\end{proposition}

From the bias - variance decomposition in the anisotropic case (see Proposition 1 in \cite{Decomp}) we get the following bound
$$\mathbb{E} [|\hat{\mu}_{h,T}(x) - \mu (x)|^2] \underset{\sim}{<} \sum_{l = 1}^d h_l^{2\beta_l} + T^{- 2}  Var ( \frac{1}{\prod_{l=1}^d h_l} \int_0^T \prod_{m = 1}^d K (\frac{x_m - X_t^m}{h_m}) dt).$$
We want to bound the variance here above using Proposition \ref{prop: variance bound} on the function \\ $f(y) := \frac{1}{\prod_{l=1}^d h_l} \prod_{m = 1}^d K (\frac{x_m - y_m}{h_m})$.
As it will be explained in the proof of Proposition \ref{prop: main result conv d>3} in Section \ref{S: Proof_Main}, for $d \ge 3$ it is
$$Var ( \frac{1}{\prod_{l=1}^d h_l} \int_0^T \prod_{m = 1}^d K (\frac{x_m - X_t^m}{h_m}) dt) \le c T(\prod_{l = 1}^d h_l)^{\frac{2}{d} - 1},  $$
which leads us to the following convergence rate.
\begin{proposition}
Suppose that A1 - A3 hold. If $\mu \in \mathcal{H}_d (\beta, \mathcal{L})$, then the estimator given in \eqref{eq: def estimator} satisfies, for $d \ge 3$, the following risk estimates:
\begin{equation}
\mathbb{E}[|\hat{\mu}_{h,T}(x) - \mu (x)|^2] \underset{\sim}{<} \sum_{l = 1}^d h_l^{2\beta_l} + T^{-1}(\prod_{l = 1}^d h_l)^{\frac{2}{d} - 1}.
\label{eq: rischio d ge 3}
\end{equation}
Defined $\frac{1}{\bar{\beta}} := \frac{1}{d} \sum_{l = 1}^d \frac{1}{\beta_l}$, the rate optimal choice $h_l = h_l (T) = (\frac{1}{T})^{\frac{\bar{\beta}}{\beta_l(2 \bar{\beta} + d -2)}}$ yields the convergence rate
$$\mathbb{E}[|\hat{\mu}_{h,T}(x) - \mu (x)|^2] \underset{\sim}{<} T^{- \frac{2\bar{\beta}}{2\bar{\beta}+ d - 2}}.$$
\label{prop: main result conv d>3} 
\end{proposition}
We underline that, in the continuous case, the convergence rate found by Strauch in \cite{Strauch} for the estimation of the invariant density $\mu$ belonging to the anisotropic H{\"o}lder class $\mathcal{H}_d (\beta + 1, \mathcal{L})$ is $T^{- \frac{2(\overline{\beta + 1})}{2(\overline{\beta + 1})+ d - 2}}$, for $d \ge 3$. \\
In Proposition \ref{prop: main result conv d>3} we estimate $\mu$ over anisotropic H{\"o}lder class $\mathcal{H}_d (\beta, \mathcal{L})$ and we therefore extend \cite{Strauch} to the jumps - diffusion case: the convergence rate we obtain is the same it was in the case without jumps, which is also analogous to the rate first obtained by Reiss and Dalalyan in \cite{RD} for the estimation of the invariant density $\mu$ over isotropic H{\"o}lder class $\mathcal{H}_d (\beta + 1, \mathcal{L})$, up to replacing the mean smoothness $\overline{\beta + 1}$ with $\beta + 1$, the common smoothness over the $d$ different dimensions. \\
 \\
For $d =1$ and $d=2$, the bound on the variance changes. Therefore, the rate optimal choice $h$ will be different as well, as explained in following two propositions.
\begin{proposition}
Suppose that A1 - A3 hold. If $\mu \in \mathcal{H}_d (\beta, \mathcal{L})$, then the estimator given in \eqref{eq: def estimator} satisfies, for $d = 1$, the following risk estimates:
\begin{equation}
\mathbb{E} [|\hat{\mu}_{h,T}(x) - \mu(x)|^2] \underset{\sim}{<} h^{2 \beta} + \frac{1}{T}(1 + (\log(\frac{1}{h}))^{2 - \frac{(1 + \alpha)}{2}} + \log(\frac{1}{h})).
\label{eq: rischio d = 1}
\end{equation}
The rate optimal choice for $h$ yields to the convergence rate
$$\mathbb{E} [|\hat{\mu}_{h,T}(x) - \mu (x)|^2] \underset{\sim}{<} \frac{(\log T)^{(2 - \frac{(1 + \alpha)}{2}) \lor 1}}{T}.$$
\label{prop: main result conv d=1} 
\end{proposition}
It is worth remarking that, in Proposition \ref{prop: main result conv d=1}, it is stated the main difference between the case with and without jumps. Indeed, if in the continuous case the convergence rate was $\frac{1}{T}$, now it depends on the degree of the jumps $\alpha$ and it is between $\frac{\log T}{T}$ and $\frac{(\log T)^\frac{3}{2}}{T}$. \\
We need to say that the convergence rate we have found here above for the estimation of the invariant density of a stochastic differential equations with jumps in the one dimensional setting is not necessarily the optimal one. In the continuous case other methods have been explored for such an estimation when $d=1$, as the use of diffusion local time to get the optimal rate $\frac{1}{T}$. We do not rule out the possibility to get a sharper bound through the exploitation of other approaches also for the jumps case, finding therefore a convergence rate faster than the one presented in the previous proposition.
\begin{proposition}
Suppose that A1 - A3 hold. If $\mu \in \mathcal{H}_d (\beta, \mathcal{L})$, then the estimator given in \eqref{eq: def estimator} satisfies, for $d = 2$, the following risk estimates:
\begin{equation}
\mathbb{E} [|\hat{\mu}_{h,T}(x) - \mu (x)|^2] \underset{\sim}{<} h_1^{2 \beta_1} + h_2^{2 \beta_2} + \frac{1}{T} (1 + \log(\frac{1}{h_1 h_2})).
\label{eq: rischio d = 2}
\end{equation}
The rate optimal choice for $h$ yields to the convergence rate
$$\mathbb{E} [|\hat{\mu}_{h,T}(x) - \mu (x)|^2] \underset{\sim}{<}\frac{\log T}{T}.$$
\label{prop: main result conv d=2} 
\end{proposition}
Comparing our result with the convergence rate obtained in the continuous case over isotropic H{\"o}lder class $\mathcal{H}_d (\beta + 1, \mathcal{L})$ in \cite{RD} and anisotropic H{\"o}lder class $\mathcal{H}_d (\beta + 1, \mathcal{L})$ in \cite{Strauch}, which is $\frac{(\log T)^4}{T}$ in both works, one can observe that the convergence rate seems being faster in presence of jumps. \\
The reason why it happens is that in \cite{RD} they assume the transition density to be upper bounded by $C (t^{- \frac{d}{2}} + t^{\frac{3d}{2}})$, which is a bound different from the one we get from Lemma \ref{lemma: bound transition}. \\
If the term $t^{\frac{3d}{2}}$ would have been absent in their assumption, e. g. for bounded drift, then the convergence rate in the continuous case could have been improved to $\frac{\log T}{T}$, which is also what we get in the jump- diffusion case. \\
In \cite{Strauch}, Nash and Poincar\'e inequalities lead the author to an upper bound on the transition density which is analogous to the one found in \cite{RD} (see Remark 2.4 of \cite{Strauch}). \\
\\
From the pointwise estimation of the invariant density gathered in the three previous propositions we move to the estimation on $L^2(A)$, where $A$ is a compact set of $\mathbb{R}^d$. \\
In the sequel, for $A \subset \mathbb{R}^d$ compact and for $g \in L^2(A)$, $\left \| g \right \|^2_A := \int_A |g(x)|^2 dx$ denotes the $L^2$ norm with respect to Lebesgue on $A$. \\
As a consequence of Propositions \ref{prop: main result conv d>3}, \ref{prop: main result conv d=1} and \ref{prop: main result conv d=2} and the fact that the constants which turn out in the proofs do not depend on $x$, the following corollary holds true:
\begin{corollary}
If $\mu \in \mathcal{H}_d (\beta, \mathcal{L})$, then for the rate optimal choice for h = h(T) provided in Propositions \ref{prop: main result conv d>3}, \ref{prop: main result conv d=1} and \ref{prop: main result conv d=2} we have the following risk estimates:
\begin{equation}
\mathbb{E}[\left \| \hat{\mu}_{h,T} - \mu \right \|^2_A] \underset{\sim}{<} V_d(T) := 
\begin{cases}
\frac{(\log T)^{(2 - \frac{(1 + \alpha)}{2}) \lor 1}}{T} \qquad \mbox{for } d=1 \\
\frac{\log T}{T} \qquad \mbox{for } d=2 \\
T^{- \frac{2\bar{\beta}}{2\bar{\beta}+ d - 2}} \qquad \mbox{for } d \ge 3.
\end{cases}
\label{eq: rates L2}    
\end{equation}
\label{cor: rate L2}
\end{corollary}

The proof of Corollary \ref{cor: rate L2} will be given in Section \ref{S: Proof_Main}.

\subsection{Adaptive procedure}\label{S: Adaptive}
The question of density estimation belongs to the canonical framework of nonparametric statistics. \\
As detailed in Propositions \ref{prop: main result conv d=1} and \ref{prop: main result conv d=2}, both the bandwidth and the upper bound on the rate of convergence appearing on the right hand side of \eqref{eq: rischio d = 1} and \eqref{eq: rischio d = 2} do not depend on the unknown smoothness of the invariant density $\mu$ and so there is no gain in implementing a data-driven bandwidth selection procedure for density estimation in the framework of continuous observations of a one or two dimensional diffusion process with jumps. Hence, throughout the sequel we restrict to the case $d \ge 3$. \\
It is clear from the previous section that for $d\ge 3$, instead, the proposed bandwidth choice depends on the regularity of the density $\mu$, which is unknown. This is why we study a data-driven bandwidth selection device. \\
We emphasize that the $d$ selected bandwidths are different, and this anisotropy property is important in our setting: the regularity in each direction can be various. The bandwidth selection procedure has to be able to provide such different choices for $h_1$, $h_2$, ... , $h_d$. \\
To select $h$ adequately, we propose the following method, inspired from Goldenshluger and Lespki \cite{Adaptive}. \\
We define the set of candidate bandwidths $\mathcal{H}_T$ as 
\begin{equation}
\mathcal{H}_T \subset \left \{ h = (h_1, ... , h_d)^T \in (0,1]^d : \, \frac{(\log T)^{2d}}{T^\frac{d}{3}} \le \prod_{l = 1}^d h_l \le ( \frac{1}{\log T})^{\frac{3 d}{d - 2}} \right \},
\label{eq: def mathcal H}
\end{equation}
The conditions on $\prod_{l = 1}^d h_l$ we have just given are needed to use Talagrand inequality, on the basis of which we show our adaptive result. \\
We suppose moreover that the growth of $|\mathcal{H}_T|$ is at most polynomial in $T$, which is there exists $c >0$ for which $|\mathcal{H}_T | \le c T^c$. \\
An example of $\mathcal{H}_T$ is the following set of candidate bandwidths:
\begin{equation}
\mathcal{H}_T := \left \{ h = (h_1, ... , h_d)^T \in (0,1]^d : \,h_i = \frac{1}{k_i} \, \mbox{ with } k_i \in \mathbb{N}, \,  \frac{(\log T)^{2d}}{T^\frac{d}{3}} \le \prod_{l = 1}^d \frac{1}{k_l} \le ( \frac{1}{\log T})^{\frac{3 d}{d - 2}} \right \}.
\label{eq: example HT}    
\end{equation}
In correspondence of the variation of $h \in \mathcal{H}_T$, we have the following family of estimators, defined as in \eqref{eq: def estimator}
$$\mathcal{F}(\mathcal{H}_T) := \left \{ \hat{\mu}_h (x) := \frac{1}{T} \int_0^T \mathbb{K}_h (X_u - x) du: \quad x \in \mathbb{R}^d, \quad h \in \mathcal{H}_T \right \}$$
where, for $y \in \mathbb{R}^d$, it is
\begin{equation}
\mathbb{K}_h (y) := \prod_{l = 1}^d \frac{1}{h_l} \prod_{m =1}^d K (\frac{y_m}{h_m}).
\label{eq: def mathbb K h}
\end{equation}
We aim at selecting an estimator from the family $\mathcal{F}(\mathcal{H}_T)$ in a completely data-driven way, based only on the observation of the continuous trajectory of the process X solution of \eqref{eq: model}. \\
We now turn to describing the selection procedure from $\mathcal{F}(\mathcal{H}_T)$, which is based on auxiliary estimators relying on the convolution operator. According to our records, it was introduced in \cite{Lep99} as a device to circumvent the lack of ordering among a set of estimators in anisotropic case, where the increase of the variance of an estimator does not imply a decrease of its bias. \\
For any bandwidths $h = (h_1, ... , h_d)^T$, $\eta = (\eta_1, ... , \eta_d)^T$ $\in \mathcal{H}_T$ and $x \in \mathbb{R}^d$, we define
$$\mathbb{K}_h * \mathbb{K}_\eta (x) := \prod_{j = 1}^d (K_{h_j} * K_{\eta_j}) (x_j) = \prod_{j = 1}^d \int_{\mathbb{R}} K_{h_j} (u - x_j) K_{\eta_j } (u) du.$$
We moreover define the kernel estimators 
$$\hat{\mu}_{h, \eta} (x) := \frac{1}{T} \int_0^T (\mathbb{K}_h * \mathbb{K}_{\eta}) (X_u - x) du, \quad x \in \mathbb{R}^d.$$
We remark that for how we have defined the kernel estimators, since the convolution is commutative, it is $\hat{\mu}_{h, \eta} = \hat{\mu}_{ \eta, h }$. \\
The proposed selection procedure relies on comparing the differences $\hat{\mu}_{h, \eta} - \hat{\mu}_{\eta}$. \\
We define
\begin{equation}
A(h) := \sup_{\eta \in \mathcal{H}_T} (\left \| \hat{\mu}_{h, \eta} - \hat{\mu}_{\eta} \right \|^2_A - V(\eta))_+,
\label{eq: def A(h)}
\end{equation}
with $$V(h) := \frac{k}{T} \, (\prod_{l = 1}^d h_l)^{\frac{2}{d} - 1}, $$
where $k$ is a numerical constant which is large. In particular, it is sufficient to choose it bigger than the constants $2 k_0^*$ and $2 k_0$ which appear in Lemma \ref{lemma: Talagrand}. Even if $k$ is not explicit, it
 can be calibrated by simulations as done for example in Section 5 of \cite{Main adapt} through the implementation of a method inspired by Goldenshluger and Lepski \cite{Adaptive} and rewritten most recently by Lacour, Massart and Rivoirard in \cite{Lacour et al}. \\ Heuristically, $A(h)$ is an estimate of the squared bias and $V(h)$ of the variance bound. It is worth noticing that the penalty term $V(h)$ which is used here comes from Proposition \ref{prop: variance bound} for the function $f$ being the Kernel function. \\
Thus, the selection is done by setting 
\begin{equation}
\tilde{h} := \mbox{arg}\min_{ h \in \mathcal{H}_T} (A (h) + V(h)).
\label{eq: def h tilde}
\end{equation}
We introduce the following notation: $\mu_h := \mathbb{K}_h * \mu$, which is the function that is estimated without bias by $\hat{\mu}_h$, i.e. $\mathbb{E}[\hat{\mu}_h(x)] = \mu_h(x)$. Moreover we define $\mu_{h, \eta} := \mathbb{K}_h * \mathbb{K}_\eta * \mu $ and a bias term
$B(h) := \left \| \mu_h - \mu \right \|^2_{\tilde{A}}$, where we have denoted as $\left \| . \right \|_{\tilde{A}}$ the $L^2$ - norm on $\tilde{A}$, a compact set in $\mathbb{R}^d$ which is such that $\tilde{A} := \left \{  \zeta \in \mathbb{R}^d : d(\zeta, A) \le 2 \sqrt{d} \right \}$. \\
The following result holds.
\begin{theorem}
Suppose that assumptions A1 - A3 hold. Then, we have
$$\mathbb{E}[\left \| \hat{\mu}_{\tilde{h}} -  \mu \right \|^2_A] \le c_1 \inf_{h \in \mathcal{H}_T} (B(h) + V(h)) + c_1 e^{ - c_2 (\log T)^2},$$
for $c_1$ and $c_2$ positive constants.
\label{th: adaptive}
\end{theorem}
The bound stated in Theorem \ref{th: adaptive} shows that the estimator leads to an automatic trade-off between the bias $\left \| \mu_h - \mu \right \|^2_{\tilde{A}}$ and the variance $V(h)$, up to a multiplicative constant $c_1$. The last term is indeed negligible.
The proof of Theorem \ref{th: adaptive} is postponed to Section \ref{S: proof adaptive}. \\
\\
We recall that Proposition \ref{prop: main result conv d>3} provides us the rate optimal choice $h(T)$ for $d \ge 3$, which is $h_l(T) = (\frac{1}{T})^{\frac{\bar{\beta}}{\beta_l ( 2 \bar{\beta} + d - 2)}}$. \\
Using such a bandwidth we will prove in Section \ref{S: proof adaptive} the following theorem.
\begin{theorem}
Suppose that assumptions A1 - A3 hold and let $\mathcal{H}_T$ be defined by \eqref{eq: example HT}. Then, we have
$$\mathbb{E}[\left \| \hat{\mu}_{\tilde{h}} -  \mu \right \|^2_A] \le c_1 (\frac{1}{T})^{\frac{2 \bar{\beta}}{ 2 \bar{\beta} + d - 2}} + c_1 e^{ - c_2 (\log T)^2},$$
for $c_1$ and $c_2$ positive constants.
\label{th: adaptive optimal}
\end{theorem}
Underlining once again that the second term in the right hand side of the equation here above is negligible compared to the first, we have that the risk estimates we get using the bandwidth provided by our selection procedure converges to zero fast. In particular, its convergence rate coincides to the optimal one provided by both \cite{RD} and \cite{Strauch} in the case without jumps.

\section{Proof convergence rates for invariant density estimation}\label{S: Proof_Main}
In this section we prove Propositions \ref{prop: main result conv d>3}, \ref{prop: main result conv d=1} and \ref{prop: main result conv d=2}, which gives us the convergence rate for the estimation of the invariant density $\mu \in \mathcal{H}_d(\beta, \mathcal{L})$ in the three different situation: $d=1$, $d=2$ and $d \ge 3$. \\
We emphasize that all the constants will appear in the proofs do not depend on the point $x$ considered. \\
We start showing the bound on the variance gathered in Proposition \ref{prop: variance bound}.
\subsection{Proof of Proposition \ref{prop: variance bound}}
\begin{proof}
We consider first of all the case $d \ge 3$.
We define the function $f_c := f - \mu (f)$. From the symmetry and the stationarity we have 
$$Var(\int_0^T f(X_s) ds) = 2 \int_0^T \int_0^s \mathbb{E} [f_c(X_s)f_c(X_t)]dt ds = 2 \int_0^T \int_0^s \mathbb{E} [f_c(X_0)f_c(X_{s -t})]dt ds $$
Applying the change if variable $u := s- t$, using Fubini and computing the integral we have that the quantity here above is equal to $2 \int_0^T (T  - u) \mathbb{E} [f_c(X_0)f_c(X_u)] du$.
Let now $0 < \delta < D \le T$, where the specific choice of $\delta$ and $D$ will be given later. The idea is to deal with the integral here above in different way for $u$ which is in different intervals. For this reason we see $\int_0^T (T  - u) \mathbb{E} [f_c(X_0)f_c(X_u)] du$ as 
\begin{equation}
\int_0^{\delta} (T  - u) \mathbb{E} [f_c(X_0)f_c(X_u)] du + \int_{\delta}^D (T  - u) \mathbb{E} [f_c(X_0)f_c(X_u)] du + \int_D^T (T  - u) \mathbb{E} [f_c(X_0)f_c(X_u)] du.
\label{eq: spezzo integrale varianza}
\end{equation}
We now observe that
\begin{equation}
\int_0^\delta (T  - u) \mathbb{E} [f_c(X_0)f_c(X_u)] du = \int_0^\delta (T  - u) (\mathbb{E} [f(X_0)f(X_u)] - (\mu (f))^2) du \le cT \int_0^\delta | < P_{u} f , f >_{\mu} | du,
\label{eq: inizio rate}
\end{equation}
where we have denoted as $< ., . >_{\mu}$ the scalar product deriving by the norm with respect to the measure $\mu$, for which $< g, h >_{\mu}:= \int_{\mathbb{R}^d}g(x) h (x) \mu (x) dx $, each $g, h \in L^2(\mu)$. In the last inequality we have moreover used that $(\mu (f))^2$ is always more than $0$. \\
Now we use Cauchy-Schwartz inequality and the fact that $P_{u} f$ is a contraction map in $L^2(\mu)$ to get
\begin{equation}
 \int_0^{\delta}| < P_{u} f , f >_{\mu} | du \le \int_0^{\delta} \sqrt{\left \| P_{u} f \right \|^2_{\mu} \left \| f \right \|^2_{\mu}} du \le  \int_0^{\delta} \sqrt{\left \| f \right \|^4_{\mu}} du \le \left \| f \right \|^2_{\infty} \mu (\mathcal{S}) \delta ,
\label{eq: conseq contraction}
\end{equation}
where in the last inequality we have used the estimation 
$$\left \| f \right \|^2_{\mu} = \int_\mathcal{S} |f(x)|^2 \mu(x) dx \le \left \| f \right \|^2_{\infty} \mu (\mathcal{S}). $$
Concerning the second integral in \eqref{eq: spezzo integrale varianza}, we remark that \eqref{eq: inizio rate} still holds on $[\delta, D]$. We then estimate it through the definition of transition semigroup. It is
\begin{equation}
\int_{\delta}^D | < P_{u} f , f >_{\mu} | du \le \int_{\delta}^D \int_{\mathbb{R}^d} |f(x)| \int_{\mathbb{R}^d} |f(y)| p_{u}(x,y) dy  \mu (x) dx du. 
\label{eq: stima density}
\end{equation}
We want to use the bound on the transition density given in Lemma \ref{lemma: bound transition} which holds for $t \in [0, T]$ but it is not uniform in $t$ big. Nevertheless, for $t \ge 1$, we have
$$p_t(x,y) = \int_{\mathbb{R}^d} p_{t - \frac{1}{2}} (x, \zeta) p_\frac{1}{2} (\zeta,y) d\zeta \le c \int_{\mathbb{R}^d} p_{t - \frac{1}{2}} (x, \zeta)(e^{- \lambda_0 (y - \zeta)^2 \frac{1}{2}} + \frac{1}{(\sqrt{\frac{1}{2}} +|y - \zeta|)^{d + \alpha}}) d\zeta \le $$
$$ \le c \int_{\mathbb{R}^d} p_{t - \frac{1}{2}} (x, \zeta) d\zeta \le c, $$
where the constant $c$ changes from line to line.
The right hand side of \eqref{eq: stima density} is therefore upper bounded by
$$ \int_{\delta}^D \int_{\mathbb{R}^d} |f(x)| c \int_{\mathbb{R}^d} |f(y)| (u^{- \frac{d}{2}} e^{- \lambda_0 \frac{|y - x|^2}{u}} + \frac{u}{(u^\frac{1}{2} + |y - x|)^{d + \alpha}} + 1) dy \, \mu (x) dx \, du \le $$
$$\le \int_{\delta}^D \int_{\mathcal{S}} |f(x)| c \int_{\mathcal{S}} |f(y)| (u^{- \frac{d}{2}} + u^{1 - \frac{(d + \alpha)}{2}} + 1) dy \, \mu (x) dx \, du \le c \left \| f \right \|^2_{\infty} \mu (\mathcal{S})|\mathcal{S}| \int_{\delta}^D (u^{- \frac{d}{2}} + u^{1 - \frac{(d + \alpha)}{2}} + 1) du,$$
where we have bounded in both integrals the absolute value of $f$ with its infinity norm. \\
Now we want to calculate the integral with respect to the variable $u$. We observe that, since $d \ge 3$, $1 - \frac{d}{2}<0$. The exponent of the second term in the integral here above, after having integrated, is $2 - \frac{d + \alpha}{2}$. It is more than zero if $d < 4 - \alpha$, which is possible only if $\alpha \in (0,1)$ and $d =3$, less then zero otherwise. \\
Therefore, we have to consider the two different possibilities, according to the fact that the exponent would be positive or negative. It follows
\begin{equation}
\int_{\delta}^D | < P_{u} f , f >_{\mu} | du  \le c \left \| f \right \|^2_{\infty} \mu_b (\mathcal{S})|\mathcal{S}|(\delta^{1 - \frac{d}{2}} + \delta^{2 - \frac{d + \alpha}{2}} 1_{\left \{ d \ge 4 - \alpha \right \}} + D^{2 - \frac{d + \alpha}{2}} 1_{\left \{ d < 4 - \alpha \right \}} + D).
\label{eq: conseq bound density}
\end{equation}
We are now left to estimate the third integral of \eqref{eq: spezzo integrale varianza}. From Lemma \ref{lemma: beta mixing} it follows it is upper bounded by
\begin{equation}
c T \int_{D}^T |< P_{u} f_c , f_c >_{\mu} | du \le c T \int_{D}^T \sqrt{\left \| P_{u} f_c \right \|^2_{L^1(\mu)} \left \| f_c \right \|^2_{\infty}} du \le c \, T  \int_{D}^T \sqrt{(e^{- \rho u} \left \| f_c \right \|_\infty )^2 \left \| f_c \right \|^2_{\infty}} du.
\label{eq: int D T}
\end{equation}
We recall it is $f_c (x) = f(x) - \mu (f)$, where $\mu (f) = \int_{\mathcal{S}} f(x) \mu (x) dx $. \\
Therefore we have 
\begin{equation}
|f_c (x)| \le |f (x)| + |\mu(f)| \le |f (x)| + \left \| f \right \|_\infty \mu(\mathcal{S})
\label{eq: stima fc}
\end{equation}
and so 
$$\left \| f_c \right \|_\infty \le \left \| f \right \|_\infty + \left \| f \right \|_\infty \mu(\mathcal{S}) \le c \left \| f \right \|_\infty, $$
where in the last inequality we have used the following estimation 
\begin{equation}
|\mu (\mathcal{S})| \le \left \| \mu \right \|_\infty |\mathcal{S}| \le c |\mathcal{S}|
\label{eq: estim mu(S)}
\end{equation} 
and the fact that $|\mathcal{S}| < 1$. \\
Therefore we get
\begin{equation}
c T \int_{D}^T |< P_{u} f_c , f_c >_{\mu} | du \le c \, T \left \| f \right \|_\infty^2 \int_{D}^T e^{- \rho u} du \le c \, T \left \| f \right \|_\infty^2 e^{- \rho D}.
\label{eq: conseq beta mixing}
\end{equation}
Replacing \eqref{eq: conseq contraction}, \eqref{eq: conseq bound density} and \eqref{eq: conseq beta mixing} in \eqref{eq: spezzo integrale varianza} we have that
\begin{equation}
 |\int_0^T (T  - u) \mathbb{E} [f_c(X_0)f_c(X_u)] du| \le c \, T \left \| f \right \|^2_\infty |\mathcal{S}|(\delta + |\mathcal{S}| \delta^{1 - \frac{d}{2}} +  |\mathcal{S}| \delta^{2 - \frac{d + \alpha}{2}} 1_{\left \{ d \ge 4 - \alpha \right \}} + |\mathcal{S}| D^{2 - \frac{d + \alpha}{2}} 1_{\left \{ d < 4 - \alpha \right \}} +  
\label{eq: da bilanciare}
\end{equation}
$$ + |\mathcal{S}| D) + c \, T \left \| f \right \|^2_\infty  e^{- \rho D}.  $$
We now want to choose $\delta$ and $D$ for which the estimation here above is as sharp as possible. Recalling that the exponent on $\delta $ are less than zero in the second and the third terms of the right hand side of \eqref{eq: da bilanciare}, we have that for a small choice of $\delta$ would correspond the smallness of the first term while the second and the third would be big, the opposite would hold for a big $\delta$. In the same way, the behaviour of the last two terms of the right hand side of \eqref{eq: da bilanciare} relies on the choice $D$. Aiming at balancing them, we define $\delta := |\mathcal{S}|^{\frac{2}{d}}$ and $D:= [\max (- \frac{2}{\rho} \log (|\mathcal{S}|), 1)] \land T$. Replacing them in \eqref{eq: da bilanciare} if $T > (- \frac{2}{\rho} \log (|\mathcal{S}|))$ we get
$$|\int_0^T (T  - u) \mathbb{E} [f_c(X_0)f_c(X_u)] du| \le c \, T \left \| f \right \|^2_\infty (|\mathcal{S}|^{1 + \frac{2}{d}} + |\mathcal{S}|^{1 + \frac{4 - \alpha}{d}} +|\mathcal{S}|^2(log|\mathcal{S}| )^{2 - \frac{d + \alpha }{2}} + |\mathcal{S}|^2 log|\mathcal{S}| + |\mathcal{S}|^2 ), $$
which give us the result we wanted remarking that both $2$ and $1 + \frac{4 - \alpha}{d}$ are always more than $1 + \frac{2}{d}$ for $d \ge 3$ and $\alpha \in (0,2)$. \\
Otherwise, if $T \le(- \frac{2}{\rho} \log (|\mathcal{S}|))$, by the definition of $D$ we obtain $D = T$.
We still have $|\mathcal{S}|^2 D^{2 - \frac{d + \alpha}{2}} 1_{\left \{ d < 4 - \alpha \right \}} \le c |\mathcal{S}|^2 (log|\mathcal{S}| )^{2 - \frac{d + \alpha }{2}}$ and, moreover, the last integral which we dealt with in \eqref{eq: conseq beta mixing} is in this case between $T$ and $T$ and so its contribution is null. Hence, the result still holds true. \\
\\
We now consider the case $d=1$. We can act exactly like we did in the case $d \ge 3$, splitting the integral in three parts. Estimations \eqref{eq: conseq contraction} and \eqref{eq: conseq beta mixing} still holds while, using again Lemma \ref{lemma: bound transition} on the interval $[\delta, D]$, we have 
$$|\int_{\delta}^D (T  - u) \mathbb{E} [f_c(X_0)f_c(X_u)] du| \le c T  \left \| f \right \|^2_{\infty} \mu (\mathcal{S})|\mathcal{S}| \int_{\delta}^D (u^{- \frac{1}{2}} + u^{1 - \frac{(1 + \alpha)}{2}} + 1) du \le $$
$$ \le c T \left \| f \right \|^2_{\infty} \mu (\mathcal{S})|\mathcal{S}|(D^\frac{1}{2} + D^{2 - \frac{(1 + \alpha)}{2}} + D),$$
where we have used that now, integrating, both the exponent we get are positive.
In total in the case $d = 1$, using also \eqref{eq: estim mu(S)}, we therefore have
$$Var(\int_0^T f(X_s) ds) \le c T\left \| f \right \|^2_{\infty} (|\mathcal{S}| \delta + |\mathcal{S}|^2(D + D^{2 - \frac{(1 + \alpha)}{2}}) + e^{- \rho D}).$$
As we have already done, we want to make the estimation here above as sharp as possible. This time there isn't any constraint on the smallness of $\delta$ and so we can choose directly $\delta := 0$. Regarding $D$ we observe that, if $\alpha > 1$, then $2 - \frac{(1 + \alpha)}{2} < 1$; otherwise $2 - \frac{(1 + \alpha)}{2} > 1$. In each case we have the same trade off we had in the case $d \ge 3$ and so we keep taking $D:= [\max (- \frac{2}{\rho} \log (|\mathcal{S}|), 1)] \land T$; it follows
$Var(\int_0^T f(X_s) ds) \le c T\left \| f \right \|^2_{\infty}|\mathcal{S}|^2(1 + (\log|\mathcal{S}|)^{2 - \frac{(1 + \alpha)}{2}} + \log (\frac{1}{|\mathcal{S}|})). $ \\
\\
In the case $d =2$ estimations \eqref{eq: conseq contraction} and \eqref{eq: conseq beta mixing} keep holding true. The bound on the transition density gathered in Lemma \ref{lemma: bound transition} for $d =2$ yields  
$$|\int_{\delta}^D (T  - u) \mathbb{E} [f_c(X_0)f_c(X_u)] du| \le c T \left \| f \right \|^2_{\infty} \mu (\mathcal{S})|\mathcal{S}| \int_{\delta}^D (u^{-1} + u^{1 - \frac{(2 + \alpha)}{2}} + 1) du \le$$
$$ \le c T \left \| f \right \|^2_{\infty} \mu (\mathcal{S})|\mathcal{S}| (\log(\frac{D}{\delta}) + D^{2 - \frac{(2 + \alpha)}{2}} + D), $$
having remarked that $2 - \frac{(2 + \alpha)}{2} = 1 - \frac{\alpha}{2} > 0$ because $\alpha \in (0,2)$. This entails, using also \eqref{eq: estim mu(S)}, 
$$Var(\int_0^T f(X_s) ds) \le c T\left \| f \right \|^2_{\infty} (|\mathcal{S}| \delta + |\mathcal{S}|^2(\log(\frac{D}{\delta})+ D^{2 - \frac{(2 + \alpha)}{2}} + D) + e^{- \rho D}).$$
Aiming at balancing the terms, we choose again $\delta := |\mathcal{S}|$ and $D:= [\max (- \frac{2}{\rho} \log (|\mathcal{S}|), 1)] \land T$. \\
It follows that $\log(\frac{D}{\delta}) \le c | \log|\log |\mathcal{S}||| +c |\log(\frac{1}{|\mathcal{S}|})|$ and so, observing that $\log|\log (|\mathcal{S}|)|$ is negligible compared to $\log(\frac{1}{|\mathcal{S}|})$, the bound on the variance becomes
$$Var(\int_0^T f(X_s) ds) \le c T\left \| f \right \|^2_{\infty}|\mathcal{S}|^2 (1 + \log(\frac{1}{|\mathcal{S}|})),$$
where we have also used that $2 - \frac{(2 + \alpha)}{2}$ is always less than $1$ and so $(\log(\frac{1}{|\mathcal{S}|}))^{2 - \frac{(2 + \alpha)}{2}} < \log(\frac{1}{|\mathcal{S}|})$. \\
The proposition is therefore proved.
\end{proof}

\subsection{Proof of Proposition \ref{prop: main result conv d>3}}
\begin{proof}
Estimation \eqref{eq: rischio d ge 3} is a straightforward consequence of the bias - variance decomposition and Proposition \ref{prop: variance bound} applied to $f(y) := \frac{1}{\prod_{l=1}^d h_l} \prod_{m = 1}^d K (\frac{x_m - y_m}{h_m})$, whose support $\mathcal{S}$ is such that $|\mathcal{S}| \le c \prod_{l=1}^d h_l$ and which is by construction such that $\left \| f \right \|_{\infty} \le c (\prod_{l=1}^d h_l)^{-1}$. \\
To find the optimal choice of $h$ we define $h_l (T) := (\frac{1}{T})^{a_l}$ for $l \in \left \{ 1, ... , d \right \}$ and we look for $a_1$, ... $a_d$ such that the upper bound of the mean-squared error in the right hand side of \eqref{eq: rischio d ge 3} would result as small as possible. \\
Replacing the definition of $h_l (T)$ in the bias - variance decomposition, it means searching for $a_1, ... , a_d$ for which we get the balance and so we have to resolve the following system:
$$
\begin{cases}
\beta_i a_i = \beta_{i + 1} a_{i + 1} \qquad \forall i \in \left \{ 1, ... , d-1 \right \} \\
2 \beta_d a_d = 1 + (\frac{2}{d} - 1) \sum_{l = 1}^d a_l.
\end{cases}   
$$
We observe that, as a consequence of the first $d-1$ equations, we can write $a_l$ as $\frac{\beta_d}{\beta_l} a_d$ for each $l \in \left \{ 1, ... , d-1 \right \} $. Therefore, the last equation becomes
$$2 \beta_d a_d = 1 + (\frac{2}{d} - 1) \beta_d a_d \sum_{l = 1}^d  \frac{1}{\beta_l}.$$
Defining $\displaystyle{\frac{1}{\bar{\beta}} := \frac{1}{d} \sum_{l = 1}^d  \frac{1}{\beta_l}}$, it follows $\displaystyle{2 \beta_d a_d = 1 + (\frac{2}{d} - 1) \beta_d a_d  \frac{d}{\bar{\beta}}}$, which yields
$$a_d = \frac{\bar{\beta}}{\beta_d (2 \bar{\beta} + (d - 2))}$$
and, replacing it in the system, we have
$$a_l = \frac{\bar{\beta}}{\beta_l (2 \bar{\beta} + (d - 2))} \quad \forall l \in \left \{ 1, ... , d \right \}.$$
Taking in the right hand side of \eqref{eq: rischio d ge 3} the rate optimal choice $h_l (T)= (\frac{1}{T})^{\frac{\bar{\beta}}{\beta_l (2 \bar{\beta} + (d - 2))}}$ we get the convergence rate wanted.
\end{proof}

\subsection{Proof of Proposition \ref{prop: main result conv d=1}}
\begin{proof}
The upper bound of the mean-squared error follows from the decomposition bias - variance and from Proposition \ref{prop: variance bound}, recalling that for $f (X_t) := \frac{1}{h} K(\frac{x - X_t}{h})$ we have $\left \| f \right \|_{\infty} \le c h^{-1}$ and its support $\mathcal{S}$ is such that $|\mathcal{S}| \le c h$. \\
Now, aiming at balancing the terms, we take $h := (\frac{1}{T})^a$; getting the mean-squared error is upper bounded by 
$$(\frac{1}{T})^{2 a \beta} + \frac{1}{T} + \frac{(a \log T)^{2 - \frac{(1 + \alpha)}{2}}}{T} + \frac{a \log T}{T}.$$
If $a$ gets bigger clearly $h$ gets smaller; it is enough to take $a$ such that $2 a \beta > 1$ to obtain the first two terms here above are negligible compared to the last ones, which gives us the convergence rate $\frac{(\log T)^{2 - \frac{(1 + \alpha)}{2}}}{T}$ for $\alpha \le 1$, $\frac{\log T}{T}$ for $\alpha > 1$.
\end{proof}

\subsection{Proof of Proposition \ref{prop: main result conv d=2}}
\begin{proof}
Again, \eqref{eq: rischio d = 2} follows naturally from the bias - variance decomposition and Proposition \ref{prop: variance bound}. \\
Regarding the convergence rate, we take again $h_l := (\frac{1}{T})^{a_l}$ for $l=1,2$. \\
It follows $\log(\frac{1}{h_1 h_2}) = a_1 \log T + a_2 \log T$ and so the mean-squared error is upper bounded by
$$(\frac{1}{T})^{2 a_1 \beta_1} + (\frac{1}{T})^{2 a_2 \beta_2} + \frac{c \log T}{T}.$$
Taking $a_1$ and $a_2$ big enough to make the first two terms here above negligible compared to the third, we get the convergence rate $\frac{\log T}{T}$. 
\end{proof}

\subsection{Proof of Corollary \ref{cor: rate L2}}
\begin{proof}
It is a straightforward consequence of Propositions \ref{prop: main result conv d>3}, \ref{prop: main result conv d=1} and \ref{prop: main result conv d=2} after having remarked that the constants which turn out in all the previous propositions do not depend on the point $x$ considered. Indeed, such propositions yield
$$\mathbb{E}[\left \| \hat{\mu}_{h,T} - \mu \right \|^2_A] = \mathbb{E}[\int_A | \hat{\mu}_{h,T} (x) - \mu (x)|^2 dx] \le c \int_A c V_d (T) dx \le c |A| V_d (T).$$
\end{proof}

\section{Proof of the adaptive procedure}\label{S: proof adaptive}
The heart of the proof of Theorem \ref{th: adaptive} consist of finding an upper bound for the expected value of $A(h)$, which is gathered in the following proposition.
\begin{proposition}
Suppose that assumptions A1 - A3 hold. Then, $\forall h \in \mathcal{H}_T$, 
$$\mathbb{E}[A (h)] \le  c_1 B(h) + c_1 e^{ - c_2 (\log T)^2}.$$
\label{prop: E[A(h)]}
\end{proposition}
Proposition \ref{prop: E[A(h)]} will be proven after the proofs of Theorems \ref{th: adaptive} and \ref{th: adaptive optimal}. 
\subsection{Proof of Theorem \ref{th: adaptive}.}
\begin{proof}
 From triangular inequality it follows $\forall h \in \mathcal{H}_T$
\begin{equation}
\left \| \hat{\mu}_{\tilde{h}} - \mu  \right \|^2 \le c( \left \| \hat{\mu}_{\tilde{h}} - \hat{\mu}_{h,\tilde{h}}  \right \|^2 + \left \| \hat{\mu}_{h,\tilde{h}} - \hat{\mu}_h  \right \|^2 + \left \| \hat{\mu}_h - \mu  \right \|^2)
\label{eq: start adapt}
\end{equation}
By the definition \eqref{eq: def A(h)} of $A(h)$ it follows that the first and the second term of \eqref{eq: start adapt} are respectively upper bounded by $A(h) + V(\tilde{h})$ and $A(\tilde{h}) + V(h)$, having also used on the second term that $\hat{\mu}_{h,\tilde{h}} = \hat{\mu}_{\tilde{h}, h}$. Then, since $\tilde{h}$ has been defined in \eqref{eq: def h tilde} as the $h \in \mathcal{H}_T$ for which $A(h) + V(h)$ is minimal, we clearly have that $A(\tilde{h}) + V(\tilde{h}) \le A(h) + V(h)$. \\
Hence, for any $h \in \mathcal{H}_T$, we get
\begin{equation}
\left \| \hat{\mu}_{\tilde{h}} - \mu  \right \|^2 \le c (A(h) + V(h) + \left \| \hat{\mu}_h - \mu  \right \|^2).
\label{eq: adapt penultima}
\end{equation}
We want an upper bound for the expected value of the left hand side of the equation \eqref{eq: adapt penultima} and so we need to evaluate $\mathbb{E} [\left \| \hat{\mu}_h - \mu  \right \|^2] = \mathbb{E}[\int_A | \hat{\mu}_h (x) - \mu(x)|^2 dx]$. \\
From the standard bias variance decomposition, recalling that $\mu_h = \mathbb{K}_h * \mu$ is such that $\mathbb{E} [\hat{\mu}_h (x)] = \mu_h (x)$, we get
$$\mathbb{E}[\int_A | \hat{\mu}_h (x) - \mu(x)|^2 dx] = \int_A | \mu_h (x) - \mu(x)|^2 dx + \int_A \mathbb{E}[| \hat{\mu}_h (x) - \mu_h(x)|^2] dx.$$
Now, we can upper bound the first term of the right hand side here above by enlarging the integration domain getting
\begin{equation}
\int_A | \mu_h (x) - \mu(x)|^2 dx \le \int_{\tilde{A}} | \mu_h (x) - \mu(x)|^2 dx = B(h).
\label{eq: estim B(h)}
\end{equation}
Moreover, as consequence of Proposition \ref{prop: variance bound} in the case $d \ge 3$ we obtain, as it was in Proposition \ref{prop: main result conv d>3},
\begin{equation}
\int_A \mathbb{E}[| \hat{\mu}_h (x) - \mu_h(x)|^2] dx = \int_A  Var( \frac{1}{T \, \prod_{l = 1}^d h_l} \int_0^T \prod_{m = 1}^d K (\frac{x_m - X_u^m }{h_m}) du ) dx \le |A| \frac{c}{T} (\prod_{l = 1}^d h_l)^{\frac{2}{d} - 1}.
\label{eq: estim adap var}
\end{equation}
Comparing the upper bound given in \eqref{eq: estim adap var} with the definition of $V(h)$ and using also \eqref{eq: adapt penultima} and \eqref{eq: estim B(h)} we get, for each $h \in \mathcal{H}_T$, 
$$\mathbb{E} [\left \| \hat{\mu}_h - \mu  \right \|^2] \le c( B(h) +  V(h)) + \mathbb{E}[A(h)].$$
Now from Proposition \ref{prop: E[A(h)]} and the arbitrariness of the bandwidth $h$ we are considering it follows
$$\mathbb{E}[\left \| \hat{\mu}_{\tilde{h}} - \mu \right \|^2] \le c_1 \inf_{h' \in \mathcal{H}_T}(B(h') + V(h')) +c_1 e^{ - c_2 (\log T)^2},$$
as we wanted.
\end{proof}

As a consequence of Theorem \ref{th: adaptive} we get, considering the rate optimal choice $h_l(T) = (\frac{1}{T})^{\frac{\bar{\beta}}{ \beta_l (2 \bar{\beta} + d - 2)}} $ provided by Proposition \ref{prop: main result conv d>3}, the estimation gathered in Theorem \ref{th: adaptive optimal}. Its proof relies on the fact that, for how we have found it in Proposition \ref{prop: main result conv d>3}, if the rate optimal bandwidth belongs to $\mathcal{H}_T$ then the $\inf_{h \in \mathcal{H}_T}(B (h) + V(h))$ is  clearly realized by it.

\subsection{Proof of Theorem \ref{th: adaptive optimal}}
\begin{proof}

We observe that, for the rate optimal choice $h(T)$ of the bandwidth, the conditions gathered in the right hand side of \eqref{eq: def mathcal H}, which are $\frac{(\log T)^{2d}}{T^\frac{d}{3}} \le \prod_{l = 1}^d h_l \le ( \frac{1}{\log T})^{\frac{3 d}{d - 2}}$,
hold true.
Indeed,
$$ \prod_{l = 1}^d h_l (T) = (\frac{1}{T})^{\frac{\bar{\beta}}{ 2 \bar{\beta} + d - 2} \sum_{l = 1}^d \frac{1}{\beta_l}} = (\frac{1}{T})^{\frac{d}{ 2 \bar{\beta} + d - 2}}. $$
the upper bound condition is therefore
$$(\frac{1}{T})^{\frac{d}{ 2 \bar{\beta} + d - 2}} \le (\frac{1}{\log T})^{\frac{3d}{d + 2}},$$
which is true if and only if $\log T \le T^{\frac{d - 2}{3(2 \bar{\beta} + d - 2)}}$. Now we observe that $\frac{d - 2}{3(2 \bar{\beta} + d - 2)} > 0$ for $\bar{\beta} > 1 - \frac{d}{2}$, which is always true since $d \ge 3$. In particular we can write $\frac{d - 2}{3(2 \bar{\beta} + d - 2)} =: \gamma \in (0,1)$ and, given that eventually for $T$ going to $\infty$ it is $\log T \le T^\gamma$, we have $\prod_{l = 1}^d h_l (T) \le (\frac{1}{\log T})^{\frac{3d}{d + 2}}$. \\
In the same way it is 
$$(\frac{1}{T})^{\frac{d}{ 2 \bar{\beta} + d - 2}} \ge \frac{(\log T)^{2 d}}{ T^{\frac{d}{3}}}.$$
For the same reasoning as here above it is true if  $(\frac{1}{3} - \frac{1}{2 \bar{\beta} + d - 2}) \frac{1}{2} =: \gamma$ is positive. \\
Being $\bar{\beta } > 1$ and $d \ge 3$, it turns out $\gamma > 0$, as we wanted. \\ 
Up to consider $\tilde{h}_l(T) := \frac{1}{\lfloor T^{\frac{\bar{\beta}}{ \beta_l (2 \bar{\beta} + d - 2)}}\rfloor }$ instead of $h_l(T)$, which is asymptotically equivalent and which leads to the same convergence rate, we have that the rate optimal choice belongs to the set of candidate bandwidths $\mathcal{H}_T$ proposed in \eqref{eq: example HT}. \\ 
Having now $h(T) \in \mathcal{H}_T$, for how we have found the rate optimal choice in Proposition \ref{prop: main result conv d>3}, the $\inf_{h \in \mathcal{H}_T} (B (h) + V (h))$ is clearly realized by it and so the bound stated in Theorem \ref{th: adaptive} is actually (see also Corollary \ref{cor: rate L2})
$$\mathbb{E}[\left \| \hat{\mu}_{\tilde{h}} -  \mu \right \|^2_A] \le c_1 (\frac{1}{T})^{\frac{2 \bar{\beta}}{ 2 \bar{\beta} + d - 2}} + c_1 e^{ - c_2 (\log T)^2},$$
as we wanted.

\end{proof}

We have showed Theorem \ref{th: adaptive} using, as main tool, the bound on $\mathbb{E} [A(h)]$ stated in Proposition \ref{prop: E[A(h)]}. Its proof, as we will see in the next section, relies on the use of Berbee's coupling method as in Viennet \cite{Viennet} and on a version of Talagrand inequality given in Klein and Rio \cite{KR}.

\subsection{Proof of Proposition \ref{prop: E[A(h)]}}
\begin{proof}
To find a bound for $\mathbb{E}[A(h)]$, for each $h \in \mathcal{H}_T$ we want to use Talagrand inequality, stated on random variables which are independent. Therefore, we start introducing some blocks which are mutually independent. We do it through the use of Berbee's coupling method as done in Viennet \cite{Viennet}, Proposition 5.1 and its proof p. 484. \\
We assume that $T = 2 p_T q_T$, with $p_T$ integer and $q_T$ a real to be chosen. We split the initial process $X = (X_t)_{t \in [0,T]}$ in $2 p_T$ processes of a length $q_T$: for each $j \in \left \{ 1, ... , p_T \right \} $ we consider \\
$X^{j,1} := (X_t)_{t \in [2(j - 1)q_T, (2j - 1)q_T]}$ and $X^{j,2} := (X_t)_{t \in [(2j-1) q_T, 2 j q_T]}$. \\
Then, there exist a process $(X^*_t)_{t \in [0,T]}$ satisfying the following properties:
\begin{enumerate}
    \item For $j \in \left \{ 1, ... , p_T \right \}$ the processes $X^{j,1}$ and $X^{* \,j,1} := (X^*_t)_{t \in [2(j - 1)q_T, (2j - 1)q_T]}$ have the same distribution and so have the processes $X^{j,2}$ and $X^{* \, j,2} := (X^*_t)_{t \in [(2j-1) q_T, 2 j q_T]}$. 
    \item For $j \in \left \{ 1, ... , p_T \right \}$, $\mathbb{P}(X^{j,1} \neq X^{* \,j,1} ) \le \beta_X(q_T)$ and $\mathbb{P}(X^{j,2} \neq X^{* \,j,2} ) \le \beta_X(q_T)$, where $\beta_X$ is the $\beta$ - mixing coefficient of the process $X$ as in \eqref{eq: def beta}.
    \item For each $k \in \left \{ 1,2 \right \}$, $X^{* \,1,k}, ... , X^{* \,p_T,k}$ are independent.
\end{enumerate}
We denote by $\hat{\mu}^*_h$ the estimator computed using $X^*_t$ instead of $X_t$ and we write $\hat{\mu}^*_h = \frac{1}{2}(\hat{\mu}^{* (1)}_h + \hat{\mu}^{*(2)}_h)$ to separate the part coming from $X^{* \,.,1}$ (super -index $(1)$) and those coming from $X^{* \,.,2}$ (super -index $(2)$), having $\hat{\mu}^{* (1)}_h := \frac{1}{ p_T q_T} \sum_{j = 1}^{p_T} \int_{2(j - 1) q_T}^{(2 j - 1) q_T}\mathbb{K}_h (X^*_u - x) du $. \\
In a natural way we define moreover $\hat{\mu}^{*}_{h, \eta} := \mathbb{K}_\eta * \hat{\mu}^*_h$, which can be written again as $\frac{1}{2}(\hat{\mu}^{* (1)}_{h, \eta} + \hat{\mu}^{*(2)}_{h, \eta})$, to separate the contribution of $X^{* \,.,1}$ and $X^{* \,.,2}$. \\
With this background we can evaluate $\mathbb{E}[A(h)]$.
We recall that, as defined in \eqref{eq: def A(h)}, $$A(h) := \sup_{\eta \in \mathcal{H}_T} (\left \| \hat{\mu}_{h, \eta} - \hat{\mu}_{\eta} \right \|^2 - V(\eta))_+.$$
Now we can see $\hat{\mu}_{h, \eta} - \hat{\mu}_{\eta}$ as sum of different terms which we deal with singularly:
$$\hat{\mu}_{h, \eta} - \hat{\mu}_{\eta} := (\hat{\mu}_{h, \eta} - \hat{\mu}_{h, \eta}^*) + (\hat{\mu}_{h, \eta}^* - \mu_{h, \eta}) + (\mu_{h, \eta} - \mu_\eta) + (\mu_\eta - \hat{\mu}_{\eta}^*) + (\hat{\mu}_{\eta}^* - \hat{\mu}_{\eta}). $$
As a consequence of the triangular inequality and of the definition of $A(h)$ the following estimation holds true:
$$A(h) \le \sup_{\eta \in \mathcal{H}_T} [\left \| \hat{\mu}_{h, \eta} - \hat{\mu}_{h, \eta}^* \right \|^2  + (\left \| \hat{\mu}_{h, \eta}^* - \mu_{h, \eta} \right \|^2 - \frac{V(\eta)}{2})_+ + \left \| \mu_{h, \eta} - \mu_\eta \right \|^2 + (\left \| \mu_\eta - \hat{\mu}_{\eta}^* \right \|^2 - \frac{V(\eta)}{2})_+ + \left \| \hat{\mu}_{\eta}^* - \hat{\mu}_{\eta} \right \|^2]=$$
$$= \sup_{\eta \in \mathcal{H}_T}[ \sum_{j = 1}^5 I_j^{h,\eta}]$$
We start considering $I_5^{h, \eta}$. We define the set 
$$\Omega^* := \left \{ X_t = X^*_t \quad \forall t \in [0,T] \right \}.$$
As a consequence of the second property of the process $X^*$ and of the $\beta$ - mixing with exponential decay showed in Lemma \ref{lemma: beta mixing} we get, recalling that $2 p_T q_t = T$ (with $q_T$ and $p_T$ to be chosen),
\begin{equation}
\mathbb{P}(\Omega^{*c}) \le 2 p_T \beta_X(q_T) \le c \frac{T}{q_T}e^{- \gamma q_T}.
\label{eq: P(Omega compl)}
\end{equation}
From the definition of $\hat{\mu}_{h}^*$ and Jensen inequality it is 
$$\left \| \hat{\mu}_{\eta}^* - \hat{\mu}_{\eta} \right \|^2 = \int_A (\frac{1}{T} \int_0^T \mathbb{K}_\eta(X_t - x) - \mathbb{K}_\eta(X^*_t - x) dt)^2 dx 1_{\Omega^{* c}} \le $$
\begin{equation}
\le c \int_A \frac{1}{T^2} T (\int_0^T \mathbb{K}^2_\eta(X_t - x) + \mathbb{K}^2_\eta(X^*_t - x) dt ) dx 1_{\Omega^{* c}} \le c \left \| \mathbb{K}_\eta \right \|_\infty^2 1_{\Omega^{* c}}
\label{eq: bound muh - muh*}
\end{equation}
By the definition \eqref{eq: def mathbb K h} we get that, $\forall h \in \mathcal{H}_T$, $\left \| \mathbb{K}_h \right \|_\infty \le (\prod_{l = 1}^d h_l)^{- 1}$. \\
We recall that, from how we have defined in \eqref{eq: def mathcal H} the set $\mathcal{H}_T$, we have that $\forall h \in \mathcal{H}_T$\\
$\prod_{l = 1}^d h_l > \frac{(\log T)^{2 d}}{T^{\frac{d}{3}}}$ and so
$$\left \| \mathbb{K}_\eta \right \|_\infty^2 < \frac{1}{(\prod_{l = 1}^d \eta_l)^2} \le \frac{T^{\frac{2 d}{3}}}{(\log T)^{4 d}}.$$
Replacing this bound in \eqref{eq: bound muh - muh*} it follows
$$ \sup_{\eta \in \mathcal{H}_T} \left \| \hat{\mu}_{\eta}^* - \hat{\mu}_{\eta} \right \|^2 \le c \frac{T^{\frac{2 d}{3}}}{(\log T)^{4 d}}1_{\Omega^{* c}}. $$
We take its expectation and we use \eqref{eq: P(Omega compl)}, getting a term which depends on $q_T$, a real to be chosen. From the arbitrariness of $q_T$ we get a convergence to zero as fast as we want, for $T$ going to $\infty$. Indeed, taking $q_T := (\log T)^2 $ yields $\forall h, \eta \in \mathcal{H}_T$, 
\begin{equation}
\mathbb{E}[ |\sup_{\eta \in \mathcal{H}_T} I_5^{h, \eta} |] = \mathbb{E}[ |\sup_{\eta \in \mathcal{H}_T}  \left \| \hat{\mu}_{\eta}^* - \hat{\mu}_{\eta} \right \|^2|] \le c \frac{T^{\frac{2 d}{3} + 1}}{(\log T)^{4 d + 2}} e^{- \gamma (\log T)^2}. 
\label{eq: I5}
\end{equation}
Regarding $\sup_{\eta \in \mathcal{H}_T} I_1^{h , \eta}$, we estimate it through \eqref{eq: I5} and the following lemma, which will be proven in the appendix. \\
\begin{lemma}
Let $f : \mathbb{R}^d \rightarrow \mathbb{R}$ be a bounded, measurable function with support $\mathcal{S}$ satisfying \\ $diam(\mathcal{S}) \le 2 \sqrt{d}$; $\tilde{A}$ a compact set in $\mathbb{R}^d$ such that $A \subset \tilde{A}$ and $\tilde{A} = \left \{ \zeta : d(\zeta, A) \le 2 \sqrt{d} \right \}$ and $g$ a function in $L^2(\tilde{A})$. Then,
$$\left \| f*g \right \|_A \le \left \|f \right \|_{1, \mathbb{R}^d} \left \| g \right \|_{2, \tilde{A}},$$
where we have denoted as $\left \| . \right \|_A$ the usual $L^2$ norm on $A$, $\left \| . \right \|_{1 , \mathbb{R}^d}$ the $L^1$ norm on $\mathbb{R}^d$ and $\left \| . \right \|_{2,\tilde{A}}$ the $L^2$ norm on $\tilde{A}$.
\label{lemma: convolution}
\end{lemma}
We recall that $\hat{\mu}_{h, \eta} = \mathbb{K}_\eta * \hat{\mu}_h$ and $\hat{\mu}^*_{h, \eta} = \mathbb{K}_\eta * \hat{\mu}^*_h$. Therefore, remarking that $diam(K) \le 2$ and so by the definition of $\mathbb{K}_\eta$ it is $diam(\mathbb{K}_\eta) \le 2 \sqrt{d}$, we can use Lemma \ref{lemma: convolution}, which yields
$$\sup_{\eta \in \mathcal{H}_T} I_1^{h, \eta} = \sup_{\eta \in \mathcal{H}_T} \left \| \mathbb{K}_\eta *(\hat{\mu}_h - \hat{\mu}^*_h) \right \|^2 \le  \sup_{\eta \in \mathcal{H}_T} \left \| \mathbb{K}_\eta \right \|_{1, \mathbb{R}^d}^2 \left \| \hat{\mu}_h - \hat{\mu}^*_h \right \|^2_{\tilde{A}}.$$
Taking the expected value, using that $\left \| \mathbb{K}_\eta \right \|_{1, \mathbb{R}^d} \le c $ $\forall \eta \in \mathcal{H}_T$ and the equation \eqref{eq: I5}, remarking that the dependence on the integration set considered is hidden in the constant $c$ in which this time will appear $|\tilde{A}|$ instead of $|A|$ we get
\begin{equation}
\mathbb{E}[ |\sup_{\eta \in \mathcal{H}_T} I_1^{h, \eta} |] \le c \frac{T^{\frac{2 d}{3} + 1}}{(\log T)^{4 d + 2}} e^{- \gamma (\log T)^2}.
\label{eq: I1}
\end{equation}
We still use Lemma \ref{lemma: convolution} to study $\sup_{\eta \in \mathcal{H}_T} I_3^{h, \eta}$, recalling that $\mu_{h, \eta}= \mathbb{K}_\eta * \mu_h$ and $\mu_\eta = \mathbb{K}_\eta * \mu_b$. It yields
\begin{equation}
\sup_{\eta \in \mathcal{H}_T} I_3^{h, \eta} = \sup_{\eta \in \mathcal{H}_T} \left \| \mathbb{K}_\eta *(\mu_h - \mu) \right \|^2 \le \sup_{\eta \in \mathcal{H}_T} \left \| \mathbb{K}_\eta \right \|^2_{1, \mathbb{R}^d} \left \| \mu_h - \mu \right \|^2_{\tilde{A}} \le c  \left \| \mu_h - \mu \right \|^2_{\tilde{A}} = c B(h).
\label{eq: I3}
\end{equation}
We are left to study $I_2^{h,\eta}$ and $I_4^{h,\eta}$, for which we need the following lemma that will be showed right after the proof of this proposition. \\
\begin{lemma}
For $i = 1,2$, there exist some positive constants $c_1^*$, $c_2^*$, $c_3^*$ and a constant $k_0^*$ such that, for any $\bar{k} \ge k_0^*$,
\begin{equation}
\mathbb{E} [\sup_{\eta \in \mathcal{H}_T} (\left \| \mu_\eta  - \hat{\mu}^{* (i)}_\eta \right \|^2 - \frac{\bar{k}}{T} (\prod_{l = 1}^d \eta_l)^{\frac{2}{d} - 1})_+] \le c_1^* T^{c_2^*} e^{- c_3^*(\log T)^2 }.
\label{eq: Talagrand per I4}
\end{equation}
Moreover there exist $k_0$ such that, for any $\bar{k} \ge k_0$, 
\begin{equation}
\mathbb{E}[\sup_{\eta \in \mathcal{H}_T} (\left \| \mu_{h,\eta}  - \hat{\mu}^{* (i)}_{h,\eta} \right \|^2 - \frac{\bar{k}}{T} (\prod_{l = 1}^d \eta_l)^{\frac{2}{d} - 1})_+] \le c_1^* T^{c_2^*} e^{- c_3^*(\log T)^2 }.
\label{eq: Talagrand per I2}
\end{equation}
\label{lemma: Talagrand}
\end{lemma}
Concerning $I_4^{h, \eta}$ observe it is $\mu_\eta  - \hat{\mu}^{*}_\eta = \frac{1}{2}(2 \mu_\eta - \hat{\mu}^{* (1)}_\eta - \hat{\mu}^{* (2)}_\eta  )$. Hence, from triangular inequality and the definition of positive part function, we get
$$I_4^{h, \eta} \le c(\left \| \mu_\eta - \hat{\mu}_{\eta}^{* (1)} \right \|^2 - \frac{V(\eta)}{2})_+ + c (\left \| \mu_\eta - \hat{\mu}_{\eta}^{* (2)} \right \|^2 - \frac{V(\eta)}{2})_+.$$
From \eqref{eq: Talagrand per I4}, for a $k$ in the definition of $V(\eta)$ big enough, for which we have $\frac{k}{2} > (k_0^* \lor k_0)$, we get
\begin{equation}
\mathbb{E} [\sup_{\eta \in \mathcal{H}_T} I_4^{h, \eta}] \le c_1 T^{c_2} e^{- c_3(\log T)^2 }.
\label{eq: I4}
\end{equation}
In the same way, remarking that 
$$I_2^{h, \eta} \le c(\left \| \hat{\mu}_{h, \eta}^{* (1)} - \mu_{h,\eta} \right \|^2 - \frac{V(\eta)}{2})_+ + c (\left \| \hat{\mu}_{h, \eta}^{* (2)} - \mu_{h,\eta} \right \|^2 - \frac{V(\eta)}{2})_+$$
and using \eqref{eq: Talagrand per I2} it follows
\begin{equation}
\mathbb{E} [\sup_{\eta \in \mathcal{H}_T} I_2^{h, \eta}] \le c_1 T^{c_2} e^{- c_3(\log T)^2 }.
\label{eq: I2}
\end{equation}
From \eqref{eq: I5} , \eqref{eq: I1}, \eqref{eq: I3}, \eqref{eq: I4} and \eqref{eq: I2} we obtain, for any $h \in \mathcal{H}_T$,
$$\mathbb{E}[A (h)] \le c \frac{T^{\frac{2 d}{3} + 1}}{(\log T)^{4 d + 2}} e^{- \gamma (\log T)^2} + c B(h) + c_1 T^{c_2} e^{- c_3(\log T)^2 } \le c B(h) + c_1 e^{- c_2(\log T)^2 },$$
as we wanted.
\end{proof}

To conclude the proof of the adaptive procedure we need to show Lemma \ref{lemma: Talagrand}, which core is the use of the Talagrand inequality. First of all, we recall the following version of the Talagrand inequality, which has been stated as Lemma 2 in \cite{Main adapt} and which is a straightforward consequence of the Talagarand inequality given in Klein and Rio \cite{KR}.
\begin{lemma}
Let $T_1, ... , T_n$ be independent random variables with values in some Polish space $\mathcal{X}$, $\mathcal{R}$ a countable class of measurable functions from $\mathcal{X}$ into $[-1, 1]^p$ and 
$v_p (r) := \frac{1}{p} \sum_{j = 1}^p [r (T_j) - \mathbb{E}[r(T_j)]]. $
Then,
\begin{equation}
\mathbb{E}[(\sup_{r \in \mathcal{R}} |v_p(r)|^2 - 2H^2)_+] \le c (\frac{v}{p} e^{- c \frac{p H^2}{v}} + \frac{M^2}{p^2} e^{- c \frac{p H}{M}}),
\label{eq: Talagrand in Klein Rio}
\end{equation}
with $c$ a universal constant and where
$$\sup_{r \in \mathcal{R}} \left \| r \right \|_\infty \le M, \quad \mathbb{E}_b[\sup_{r \in \mathcal{R}}|v_p (r)|] \le H, \quad \sup_{r \in \mathcal{R}} \frac{1}{p} \sum_{j = 1}^p Var(r (T_j)) \le v.$$
\label{lemma: Talagrand in Klein Rio}
\end{lemma}

\subsection{Proof of Lemma \ref{lemma: Talagrand}}
\begin{proof}
Since the two cases $i = 1$ and $i =2$ are similar, we study only one of them. We start proving \eqref{eq: Talagrand per I4}, the proof of inequality \eqref{eq: Talagrand per I2} follows the same line. 
We first observe it is 
\begin{equation}
\mathbb{E} [\sup_{\eta \in \mathcal{H}_T} (\left \| \mu_\eta  - \hat{\mu}^{* (1)}_\eta \right \|_A^2 - \frac{\bar{k}}{T} (\prod_{l = 1}^d \eta_l)^{\frac{2}{d} - 1})_+] \le \sum_{\eta \in \mathcal{H}_T} \mathbb{E} [(\left \| \mu_\eta  - \hat{\mu}^{* (1)}_\eta \right \|_A^2 - \frac{\bar{k}}{T} (\prod_{l = 1}^d \eta_l)^{\frac{2}{d} - 1})_+].
\label{eq: start Talagrand}
\end{equation}
Our goal is now to find a bound for the right hand side of the inequality here above using the version of the Talagrand inequality gathered in Lemma \ref{lemma: Talagrand in Klein Rio}. To do it, we need to introduce some notation. \\
We observe that $\left \| \mu_\eta  - \hat{\mu}^{* (1)}_\eta \right \|_A^2 = \sup_{r, \left \| r \right \| = 1} < \mu_\eta  - \hat{\mu}^{* (1)}_\eta,r>^2$, and the supremum can be considered over a countable dense set of function $r$ such that $\left \| r \right \| = 1$; let us denote this set by $\mathcal{B}(1)$. \\
We define
$$T_j (z) := \frac{1}{q_T} \int_{(2 j - 1) q_T}^{2j q_T} \mathbb{K}_\eta (X_t^{* j, 1} - z) dt; \quad r(T_j) := \int_{A} T_j (z) r (z) dz.$$
Thus
$$v_{p_T} (r) = <\hat{\mu}^{* (1)}_\eta - \mu_\eta,r> = \frac{1}{p_T} \sum_{j = 1}^{p_T} [r (T_j)- \mathbb{E} [r (T_j)]]$$
is a centered empirical process with independent variables 
$$\psi_r (X^{* j, 1}) := r (T_j)- \mathbb{E}[r (T_j)] = \int_{A} (\frac{1}{q_T} \int_{(2 j - 1) q_T}^{2j q_T} [\mathbb{K}_\eta (X_t^{* j, 1} - z) - \mathbb{E} [\mathbb{K}_\eta (X_t^{* j, 1} - z)]] dt) r(z) dz, $$
to which we want to apply Talagrand inequality \eqref{eq: Talagrand in Klein Rio}. Therefore, we have to compute $M$, $H$ and $v$ as defined in Lemma \ref{lemma: Talagrand in Klein Rio}. 
We start by the calculation of M. For any $r \in \mathcal{B}(1)$ it is, using the definition of $r$ and Cauchy - Schwartz inequality,
$$|\int_A T_j(z) r (z) dz| \le (\int_A T_j^2 (z) dz )^\frac{1}{2}.$$
Now from the definition of $T$ and Jensen inequality it follows
$$\int_A T_j^2 (z) dz \le \int_A \frac{1}{q_T^2} q_T \int_{(2 j - 1) q_T}^{2j q_T} \mathbb{K}^2_\eta (X_t^{* j, 1} - z) dt dz \le c \left \| K_\eta \right \|^2_\infty |\mathcal{S}| \le  c (\prod_{l = 1}^d \eta_l)^{-2} (\prod_{l = 1}^d \eta_l) = c(\prod_{l = 1}^d \eta_l)^{-1},$$
where we have also used that the support of $K_\eta$ is on $\mathcal{S}$ which size is $\prod_{l = 1}^d \eta_l$. Hence,
\begin{equation}
|\int_A T_j(z) r (z) dz| \le c(\prod_{l = 1}^d \eta_l)^{-\frac{1}{2}} =: M.
\label{eq: computation M}
\end{equation}
Regarding the computation of $H$, from the definition of $v_{p_T} (r)$ and the fact that the random variables $\psi_r^{* j, 1}$ are centered and independents it follows
$$\mathbb{E} [\sup_{r \in \mathcal{B}(1)}|v_{p_T}(r)|^2] = \mathbb{E}_b [\left \| \mu_\eta  - \hat{\mu}^{* (1)}_\eta \right \|_A^2] = \int_A Var(\frac{1}{p_T} \sum_{j = 1}^{p_T}\frac{1}{q_T} \int_{(2 j - 1) q_T}^{2j q_T} (\mathbb{K}_\eta (X_t^{* j, 1} - z) - \mathbb{E} [\mathbb{K}_\eta (X_t^{* j, 1} - z)]) dt) dz = $$
$$= \int_A \frac{1}{p_T} Var(\frac{1}{q_T} \int_{(2 j - 1) q_T}^{2j q_T} \mathbb{K}_\eta (X_t^{* j, 1} - z) dt) dz \le c |A| \frac{1}{p_T} \frac{1}{q_T}(\prod_{l = 1}^d \eta_l)^{\frac{2}{d} - 1}, $$
where in the last inequality we have used the estimation for the variance in the case $d \ge 3$ gathered in Proposition \ref{prop: variance bound}, considering that taking the Kernel function as $f$ we have that its support $\mathcal{S}$ is such that $|\mathcal{S}| \le c (\prod_{l = 1}^d \eta_l)$. It yields
\begin{equation}
\mathbb{E} [\sup_{r \in \mathcal{B}(1)}|v_{p_T}(r)|^2] \le \frac{c}{T} (\prod_{l = 1}^d \eta_l)^{\frac{2}{d} - 1} =: H^2 
\label{eq: computation H}
\end{equation}
In order to use Lemma \ref{lemma: Talagrand in Klein Rio}, we are left to compute $v$.\\
We observe it is 
$$\frac{1}{p_T} \sum_{j = 1}^{p_T} Var(\int_A \frac{1}{q_T} \int_{(2 j - 1) q_T}^{2j q_T} \mathbb{K}_\eta (X_t^{* j, 1} - z) dt \, r(z) dz) = \frac{1}{p_T} \sum_{j = 1}^{p_T} Var( \frac{1}{q_T} \int_{(2 j - 1) q_T}^{2j q_T} (\mathbb{K}_\eta * r) (X_t^{* j, 1}) dt).$$
We want to prove a tight upper bound for the variance of the integral functional $\frac{1}{q_T} \int_{(2 j - 1) q_T}^{2j q_T} f_{\eta} (X_t^{* j, 1}) dt$ of the diffusion $X^*$ , where we have denoted $f_\eta := \mathbb{K}_\eta * r$. \\
Following the proof we have given of Proposition \ref{prop: variance bound} we have, 
$$Var(\frac{1}{q_T} \int_{(2 j - 1) q_T}^{2j q_T} f_{\eta} (X_t^{* j, 1}) dt) \le \frac{c}{q_T^2} \int_0^{q_T} (q_T - u) \mathbb{E}[f_{\eta,c}(X_0^{* j, 1}), f_{\eta,c}(X_u^{* j, 1})] du =$$
$$ = \frac{c}{q_T^2} (\int_0^{D} (q_T - u) \mathbb{E}[f_{\eta,c}(X_0^{* j, 1}), f_{\eta,c}(X_u^{* j, 1})] du + \int_D^{q_T} (q_T - u) \mathbb{E}[f_{\eta,c}(X_0^{* j, 1}), f_{\eta,c}(X_u^{* j, 1})] du),$$
where we have introduced $f_{\eta,c}(x) := f_{\eta} (x) - \mu (f_\eta)$ and $D$ a quantity to be chosen in order to balance the contribution of the two integrals here above. We denote moreover as $P^*_{t}$ the transition semigroup of the process $X^*$ Concerning the integral between $0$ and $D$, we act like we did in \eqref{eq: inizio rate} and we use Cauchy - Schwartz inequality and the fact that $P^*_{ t} $ is a contraction. It follows
$$\frac{c}{q_T^2} \int_0^{D} (q_T - u) \mathbb{E}[f_{\eta,c}(X_0^{* j, 1}), f_{\eta,c}(X_u^{* j, 1})] du \le \frac{c}{q_T} |\int_0^{D} <P^*_{t} f_\eta, f_\eta>_{\mu} dt| \le$$
\begin{equation}
\le \frac{c}{q_T} \int_0^D (\left \| P^*_{t} f_\eta \right \|^2_{\mu} \left \| f_\eta \right \|^2_{\mu})^\frac{1}{2} dt \le \frac{c}{q_T} \int_{0}^D \left \| f_\eta \right \|^2_{\mu} dt \le \frac{cD}{q_T},
\label{eq : int 0 D}
\end{equation}
where in the last inequality we have used the fact that $\left \| \mu \right \|_\infty \le c$, Young inequality and the definition of the Kernel function and of $r$ in order to say
\begin{equation}
\left \| f_\eta \right \|^2_{\mu} = \left \| \mathbb{K}_\eta * r \right \|^2_{\mu} \le c \left \| \mathbb{K}_\eta * r \right \|^2_{2, \mathbb{R}^d} \le c \left \| \mathbb{K}_\eta \right \|_{1, \mathbb{R}^d}^2 \left \| r \right \|^2_{2, \mathbb{R}^d} \le c. 
\label{eq: stima norma 1 f eta}
\end{equation}
Regarding the integral between $D$ and $q_T$, we act like we did in \eqref{eq: int D T} using the exponential ergodicity gathered in Lemma \ref{lemma: beta mixing} to get
\begin{equation}
\frac{c}{q_T^2} \int_D^{q_T} (q_T - u) \mathbb{E}[f_{\eta,c}(X_0^{* j, 1}), f_{\eta,c}(X_u^{* j, 1})] du \le \frac{c}{q_T} |\int_D^{q_T} <P^*_{t} f_{\eta,c}, f_{\eta,c}>_{\mu} dt| \le \frac{c}{q_T} \int_D^{q_T} e^{- \rho t} \left \| f_{\eta,c} \right \|^2_{\infty} dt.
\label{ eq: int D qT}
\end{equation}
We now recall that $f_{\eta,c} (x) = f_\eta(x) - \mu (f_\eta)$ with $\mu(f_\eta) = \int_{\mathbb{R}^d} f_\eta (x) \mu(x) dx$. Hence, from Cauchy- Schwartz inequality and \eqref{eq: stima norma 1 f eta}, we get $|\mu (f_\eta)| \le c$ and therefore
\begin{equation}
|f_{\eta,c} (x)| \le |f_{\eta}(x) |+ c.
\label{eq: stima f eta c}
\end{equation}
To estimate the infinity norm of $f_{\eta, c}$ we remark it is, $\forall y \in \mathbb{R}^d$,
$$|f_{\eta}(y) | = |\mathbb{K}_\eta * r (y) | = | \int_A \mathbb{K}_\eta (y - z) r(z) dz  | \le (\int_A \mathbb{K}^2_\eta (y - z) dz)^\frac{1}{2} (\int_A r^2 (z) dz)^\frac{1}{2} = (\int_{\mathcal{S}} \mathbb{K}^2_\eta (y - z) dz)^\frac{1}{2},$$
where we have used Cauchy - Schwartz inequality and the fact that $r \in \mathcal{B} (1)$ and so its $2$ - norm is equal to $1$ by definition. Moreover, the Kernel function is different from $0$ only on its support $\mathcal{S}$, which size is $\prod_{l = 1}^d \eta_l$. Hence, recalling also that the infinite norm of $\mathbb{K}_\eta$ is upper bounded by $c(\prod_{l = 1}^d \eta_l)^{- 1} $, we get
$$|\mathbb{K}_\eta * r (y) | \le c (\left \|  \mathbb{K}_\eta \right \|^2_\infty | \mathcal{S} |)^\frac{1}{2} \le \frac{c}{\sqrt{\prod_{l = 1}^d \eta_l}}.$$
It follows 
\begin{equation}
\left \| f_{\eta,c} \right \|_{\infty} \le  \frac{c}{\sqrt{\prod_{l = 1}^d \eta_l}} + c \le \frac{c}{\sqrt{\prod_{l = 1}^d \eta_l}},
\label{eq: norma inf f eta c}
\end{equation}
given that the second term is negligible compared to the first.
Replacing \eqref{eq: norma inf f eta c} in \eqref{ eq: int D qT}, using also \eqref{eq : int 0 D}, we obtain 
$$ Var( \frac{1}{q_T} \int_{(2 j - 1) q_T}^{2j q_T}  f_{\eta} (X_t^{* j, 1}) dt) \le \frac{c}{q_T} (D + \frac{e^{- \rho D}}{\prod_{l = 1}^d \eta_l}).$$
We look for a $D$ for which the first and the second term of the right hand side of the inequality here above have the same magnitude. Therefore, we choose $D := [\max(- \frac{1}{\rho} \log (\prod_{l = 1}^d \eta_l), 1)] \land q_T$. Replacing such a value we get, if $q_T > - \frac{1}{ \rho} \log (\prod_{l = 1}^d \eta_l)$,
$$\frac{1}{p_T} \sum_{j = 1}^{p_T} Var(\int_A \frac{1}{q_T} \int_{(2 j - 1) q_T}^{2j q_T} \mathbb{K}_\eta (X_t^{* j, 1} - z) dt \, r(z) dz) \le \frac{c}{q_T} (1 + \log (\frac{1}{|\prod_{l = 1}^d \eta_l |}))$$
Otherwise, if $q_T \le - \frac{1}{ \rho} \log (\prod_{l = 1}^d \eta_l)$, by the definition of $D$ we have $D= q_T$. We still have the contribution of $\frac{c}{q_T} D$ which is in this case less than $\frac{c}{q_T}(\log (\frac{1}{|\prod_{l = 1}^d \eta_l |}))$ and, moreover, the contribution of the integral between $D$ and $q_T$ is now null since we have $D = q_T$.
Hence we have 
\begin{equation}
 v := \frac{c}{q_T} (1 + \log (\frac{1}{|\prod_{l = 1}^d \eta_l |})).
\label{eq: computation v}
\end{equation}
We use Lemma \ref{lemma: Talagrand in Klein Rio} on the right hand side of \eqref{eq: start Talagrand}, recalling that $\left \| \mu_\eta  - \hat{\mu}^{* (1)}_\eta \right \|_A^2 = \sup_{r \in \mathcal{B}(1)} < \mu_\eta  - \hat{\mu}^{* (1)}_\eta,r>^2 = \sup_{r \in \mathcal{B}(1)}|v_{p_T}(r)|^2$; with $M$, $H$ and $v$ as found in \eqref{eq: computation M}, \eqref{eq: computation H} and \eqref{eq: computation v}. It follows
$$\mathbb{E}[\sup_{\eta \in \mathcal{H}_T} (\left \| \mu_\eta  - \hat{\mu}^{* (1)}_\eta \right \|^2 - \frac{\bar{k}}{T} (\prod_{l = 1}^d \eta_l)^{\frac{2}{d} - 1})_+] \le $$
$$ \le c \sum_{\eta \in \mathcal{H}_T} \frac{(1 + \log (\frac{1}{|\prod_{l = 1}^d \eta_l|}))}{p_T q_T} e^{- c \frac{\frac{p_T}{T} (\prod_{l = 1}^d \eta_l)^{\frac{2}{d} - 1}}{\frac{1}{q_T} (1 + \log (\frac{1}{|\prod_{l = 1}^d \eta_l |})) }} + \frac{(\prod_{l = 1}^d \eta_l)^{-1}}{p_T^2 } e^{- c \frac{p_T \frac{1}{\sqrt{T}}(\prod_{l = 1}^d \eta_l)^{\frac{1}{d} - \frac{1}{2}}}{(\prod_{l = 1}^d \eta_l)^{- \frac{1}{2}}}} .$$
We recall that $2 p_T q_T = T$, where $q_T$ is chosen above equation \eqref{eq: I5} as $ (\log T)^2 $; we
can therefore upper bound the right hand side of the equation here above with
$$c \sum_{\eta \in \mathcal{H}_T} \frac{(1 + \log (\frac{1}{|\prod_{l = 1}^d \eta_l|}))}{T} e^{- \frac{c}{(\prod_{l = 1}^d \eta_l)^{ 1- \frac{2}{d}}(1 + \log (\frac{1}{|\prod_{l = 1}^d \eta_l|}))}} + \frac{(\log T)^4}{(\prod_{l = 1}^d \eta_l) T^2} e^{- c \frac{\sqrt{T}}{(\log T)^2} (\prod_{l = 1}^d \eta_l)^{\frac{1}{d}}} \le$$
$$ \le (\frac{\log T}{T} e^{- c (\log T)^2} + \frac{T^{\frac{d}{3} - 2}}{(\log T)^{2 d - 2}} e^{- c T^{\frac{1}{6}}})|\mathcal{H}_T|,$$
where in the last inequality we have used that, by the definition \eqref{eq: def mathcal H} we have given of $\mathcal{H}_T$, $\forall h \in \mathcal{H}_T$ we have
$\frac{(\log T)^{2d}}{T^\frac{d}{3}} \le \prod_{l = 1}^d h_l \le ( \frac{1}{\log T})^{\frac{3 d}{d - 2}}$ and so $(\prod_{l = 1}^d \eta_l)^{ 1- \frac{2}{d}}(1 + \log (\frac{1}{|\prod_{l = 1}^d \eta_l|})) \le c \frac{1}{(\log T)^2}$ and $\frac{\sqrt{T}}{(\log T)^2} (\prod_{l = 1}^d \eta_l)^{\frac{1}{d}} \ge c T^\frac{1}{6}$; we have therefore upper bounded each element of the sum with a quantity which does not depend on $\eta$. \\
We have moreover assumed that $|\mathcal{H}_T|$ has polynomial growth in $T$; and so there is a constant $c > 0$ such that
$$\mathbb{E} [\sup_{\eta \in \mathcal{H}_T} (\left \| \mu_\eta  - \hat{\mu}^{* (1)}_\eta \right \|^2 - \frac{\bar{k}}{T} (\prod_{l = 1}^d \eta_l)^{\frac{2}{d} - 1})_+] \le (\frac{\log T}{T} e^{- c (\log T)^2} + \frac{T^{\frac{d}{3} - 2}}{(\log T)^{2 d - 2}} e^{- c T^{\frac{1}{6}}})T^c ;$$
inequality \eqref{eq: Talagrand per I4} follows. \\
As a consequence of Lemma \ref{lemma: convolution} and the definition of Kernel function we moreover have 
$$\left \| \mu_{h,\eta}  - \hat{\mu}^{* (1)}_{h,\eta} \right \|^2 = \left \| \mathbb{K}_h* ( \mu_\eta  - \hat{\mu}^{* (1)}_\eta) \right \|^2 \le \left \| \mathbb{K}_h \right \|^2_{1, \mathbb{R}^d} \left \| \mu_\eta  - \hat{\mu}^{* (1)}_\eta \right \|^2_{2, \tilde{A}} \le c \left \| \mu_\eta  - \hat{\mu}^{* (1)}_\eta \right \|^2_{2, \tilde{A}}.$$
Using that and \eqref{eq: Talagrand per I4} we have just shown we obtain \eqref{eq: Talagrand per I2}. 

\end{proof}

\appendix
\section{Appendix}
\subsection{Proof of Lemma \ref{lemma: bound transition}}
\begin{proof}
Lemma \ref{lemma: bound transition} relies on the first point of Theorem 1.1 in \cite{Chen}, which needs the following assumptions on the coefficients $a$, $b$, $\gamma$ and on the jumps. \\
\textbf{($H^a$)} There are $c_1 > 0$ and $\beta \in (0,1)$ such that for all $x, y \in \mathbb{R}^d$, $|a (x) - a (y)| \le c_1 |x - y|^\beta$ and, for some $c_2 \ge 1$, $c_2^{-1} \mathbb{I}_{d \times d} \le a(x) \le c_2 \mathbb{I}_{d \times d}$. \\
\textbf{($H^k$)} The function $k(x,z) := |z|^{d + \alpha} F(\frac{z}{\gamma(x)})$ is bounded, measurable and, if $\alpha =1$, for any $0 < r < R < \infty$ it is
\begin{equation}
\int_{r \le z \le R} z k(x,z) |z|^{-d -1} dz =0. 
\label{eq: cond chen k}
\end{equation}
\textbf{($H^b$)} The function $b$ belongs to the Kato class $\mathbb{K}^2 $ which is, as defined in \cite{Chen},
$$\mathbb{K}^\gamma := \left \{ f : \mathbb{R}^d \rightarrow \mathbb{R} \mbox{ satisfies} \lim_{\delta \rightarrow 0} \sup_x \int_0^\delta \int_{\mathbb{R}^d} |f(x \, \underline{+} \, y )| \eta_{\gamma , \gamma -1} (s,y) dy \, ds =0 \right \}, $$
where we have denoted $f(x \, \underline{+} \, y )$ as an abbreviation for $f (x + y) + f (x - y)$ and 
$$\eta_{\alpha, \gamma } (t,x) := t^{\frac{\gamma}{2}} (|x| + t^{\frac{1}{2}})^{- d - \alpha}.$$
Through this paper we have assumed that Assumptions A1 - A3 hold. \\
From A1 it follows $H^a$ since we have asked the boundedness of $a$ and, in the case in which $x$ and $y$ are such that $|x - y| > 1$ we have that $\exists c$ such that
$|a (x) - a (y)| \le |a (x) | + |a (y)| \le 2c \le 2c |x - y|^\beta$, for each $\beta \in (0,1)$.
When $|x - y|\le 1$, instead, we have $|a (x) - a (y)| \le L|x -y| = L |x - y|^{1 - \beta } |x - y|^\beta \le L |x - y|^\beta$, as consequence of the Lipschitz continuity.  \\
Regarding $H^k$, we have that $k$ is a bounded function as a consequence of second and third points of A3. Indeed,
$$|k(x,z) | = |z|^{d + \alpha} |F(\frac{z}{\gamma(x)})| \le |z|^{d + \alpha} \frac{|\gamma(x)|}{|z|^{d + \alpha} } \le \gamma_{max} < \infty. $$
Observing moreover that \eqref{eq: cond chen k} holds true on the basis of the fourth point of A3, $H^k$ clearly follows. \\
We are left to show that $b \in \mathbb{K}^2$.
Noticing that (cf also Remark 2.6 in \cite{Chen})
$$\int_0^\delta \eta_{\alpha , \alpha - 1} (s,y) ds = \int_0^\delta \frac{s^{\frac{(\alpha - 1)}{2}}}{(|y| + s^\frac{1}{2})^{d + \alpha}} ds \le c \frac{(|y|^2 \land \delta)^\frac{(1 + \alpha)}{2}}{|y|^{d + \alpha}} = \frac{1}{|y|^{d - 1}} (1 \land \frac{\delta}{|y|^2})^{\frac{(1 + \alpha)}{2}},$$
it is enough to show that
$$M_b^2 (\delta) := \sup_{x \in \mathbb{R}^d} \int_{\mathbb{R}^d} |b (x + y)| \frac{1}{|y|^{d - 1}} (1 \land \frac{\delta}{|y|^2})^{\frac{3}{2}} dy \longrightarrow 0 \quad \mbox{as } \delta \rightarrow 0 $$
to get $b \in \mathbb{K}^2$. \\ \\
From Assumption A1 we now $b$ is upper bounded by a constant. Therefore we have
\begin{equation}
 \sup_{x \in \mathbb{R}^d} |\int_{\mathbb{R}^d} |b (x + y)| \frac{1}{|y|^{d - 1}} (1 \land \frac{\delta}{|y|^2})^{\frac{3}{2}} dy | \le c \int_{\left \{ y : |y| \le \sqrt{\delta}  \right \} } \frac{1}{|y|^{d - 1}} dy + c \int_{\left \{ y : |y|  > \sqrt{\delta}  \right \} } \frac{1}{|y|^{d - 1}}(\frac{\delta}{|y|^2})^{\frac{3}{2}} dy.
\label{eq: proof b in K2}
\end{equation}
We move to polar coordinates system, getting the right hand side of the equation \eqref{eq: proof b in K2} here above is upper bounded by 
$$\int_0^{\sqrt{\delta}} c d\rho + c\delta^{\frac{3}{2}} \int_{\sqrt{\delta}}^\infty \rho^{- 3} d\rho = c \sqrt{\delta} + c\delta^{\frac{3}{2}} \frac{1}{2 \delta} ,$$
which clearly goes to 0 for $\delta \rightarrow 0$. 
It yields that $H^b$ holds true. \\
It entails we can use Theorem 1.1 of \cite{Chen}; Lemma \ref{lemma: bound transition} follows.

\end{proof}

\subsection{Proof of Lemma \ref{lemma: beta mixing}}
\begin{proof}
The exponential ergodicity and the exponential $\beta$ - mixing of the process $X$ are showed in Proposition 3.8 and the second point of Theorem 2.2 of \cite{18 GLM}. \\
To use them we have to show that Assumptions 1, 2 and 3* stated in \cite{18 GLM} hold. \\ 
Assumption 1 of \cite{18 GLM} is a regularity condition which corresponds to our assumption A1. \\
We want to show that point $b$ of Assumption 2 of \cite{18 GLM} holds, which is the following: \\
(b) There exists a constant $\Delta > 0$ such that $X_{\Delta}$ admits a density $p_\Delta (x, y)$ with respect to the Lebesgue measure on $\mathbb{R}^d$ for every $x \in \mathbb{R}^d$, and $(x,y) \mapsto p_\Delta (x, y)$ is bounded in $y \in \mathbb{R}^d$ and in $x \in K$ for every compact $K \subset \mathbb{R}^d$. Moreover, for every $x \in \mathbb{R}^d$ and every open ball $U \subset \mathbb{R}^d$ there exists a point $z = z (x, U) \in supp(F)$ such that $\gamma(x) \cdot z \in U$. \\
We observe that the existence of a bounded density has already been proven in Lemma \ref{lemma: bound transition}. Moreover, from second and third points of A3, we know that $supp(F) = \mathbb{R}^d$ and that $\gamma$ is an invertible matrix. Hence, for every $x \in \mathbb{R}^d$ and every open ball $U \subset \mathbb{R}^d$ there exists a point $z = z (x, U) \in \mathbb{R}^d$ such that $\gamma(x) \cdot z \in U$. \\
To conclude, we have to prove that Assumption 3* holds and so we have to show the existence of a Lyapunov function. We therefore want to provide a function $f^*$ which satisfies the drift condition $A f^* \le -c_1 f^* + c_2$ for $c_1 > 0$ and $c_2 > 0$. $A$ denotes the generator of the diffusion, which is the sum of the continuous and discrete part
$$A_c f (x) := \frac{1}{2} \sum_{i,j = 1}^d a^2_{i,j} (x) \partial_{i,j}^2 f(x) + \sum_{i = 1}^d b_i(x) \partial_i f(x) \qquad \mbox{and}$$
$$A_d f(x) := \int_{\mathbb{R}^d} (f (x + \gamma(x) \cdot z) - f(x) - \gamma(x) \cdot z \cdot \nabla f(x) ) F(z) dz,$$
for every function $f : \mathbb{R}^d \rightarrow \mathbb{R}$, $f \in \mathcal{C}^2(\mathbb{R}^d)$. \\
From the fifth point of condition A3 we know there exists $\epsilon > 0$ such that $\int_{\mathbb{R}^d} |z|^2 e^{\epsilon |z|} F(z) dz \le c$. For such an $\epsilon$ we define $f^*(x) := e^{\epsilon |x|}$. We observe it is $\partial_i f^* (x) = \epsilon e^{\epsilon |x|} \frac{x_i}{|x|}$ and
\begin{equation}
\partial^2_{i,j} f^*(x) = \epsilon \frac{x_i x_j}{|x|^2} e^{\epsilon |x|} (\epsilon - \frac{1}{|x|}) + \epsilon e^{\epsilon |x|} \frac{1}{|x|} 1_{j = i}.
\label{eq: deriv seconde f*}
\end{equation}
We therefore have, using also the drift condition gathered in assumption A2, $\forall x : |x| > \tilde{\rho}$
$$|A_c f^* (x)| \le \frac{1}{2} \epsilon e^{\epsilon |x|} (\epsilon + \frac{2}{|x|})\sum_{i,j= 1}^d |a^2_{i,j} (x)| + \epsilon e^{\epsilon |x|} \frac{1}{|x|} < x, b(x) > \, \le$$
\begin{equation}
\le c \epsilon e^{\epsilon |x|} (\epsilon + \frac{1}{|x|})\sum_{i,j= 1}^d |a^2_{i,j} (x)| - \tilde{C} \epsilon e^{\epsilon |x|}. 
\label{eq: Ac f*}
\end{equation}
Concerning the discrete part of the generator, from intermediate value theorem we have 
$$|A_d f^* (x)| \le \int_{\mathbb{R}^d} \int_0^1 (\gamma (x) \cdot z)^T \cdot H^2 f_{(x + s \gamma (x) \cdot z)} \cdot (\gamma (x) \cdot z) ds dz,$$
where $H^2 f_{(x + s \gamma (x) \cdot z)}$ denotes the hessian matrix of the function $f$ computed in the point $x + s \gamma (x) \cdot z$. \\
We now split the integral in the right hand side here above, to act differently depending on whether  $|z| $ is more or less than $ \frac{ |x|}{2 \left \| \gamma \right \|_\infty }$. We therefore get 
$$|A_d f^* (x)| \le \int_{z : |z| \le \frac{ |x|}{2 \left \| \gamma \right \|_\infty } } \int_0^1 (\gamma (x) \cdot z)^T \cdot H^2 f_{(x + s \gamma (x) \cdot z)} \cdot (\gamma (x) \cdot z) ds dz + $$
$$ + \int_{z : |z| > \frac{ |x|}{2 \left \| \gamma \right \|_\infty }} \int_0^1 (\gamma (x) \cdot z)^T \cdot H^2 f_{(x + s \gamma (x) \cdot z)} \cdot (\gamma (x) \cdot z) ds dz =: I_1 + I_2.$$
Concerning $I_1$, from \eqref{eq: deriv seconde f*} it follows
$$I_1 \le c \int_{z : |z| \le \frac{ |x|}{2 \left \| \gamma \right \|_\infty } } \int_0^1 |z|^2 \left \| \gamma \right \|_\infty^2 \epsilon e^{\epsilon |x + s z \gamma(x)|} (\epsilon + \frac{1}{|x + s z \gamma(x)|}) F(z) dz ds \le  $$
\begin{equation}
\le c \epsilon e^{\epsilon |x|} \int_{z : |z| \le  \frac{ |x|}{2 \left \| \gamma \right \|_\infty } } |z|^2 e^{\epsilon |z| \left \| \gamma \right \|_\infty } (\epsilon + \frac{1}{|x| - |z| \left \| \gamma \right \|_\infty}) F(z) dz \le c \epsilon e^{\epsilon |x|}(\epsilon + \frac{1}{|x|}),
\label{eq: I1 f*}
\end{equation}
where in the last inequality we have used the fifth point of A3 and the fact that on the integral we are considering it is $|z| \le \frac{ |x|}{2 \left \| \gamma \right \|_\infty }$ and so $\frac{1}{|x| - |z| \left \| \gamma \right \|_\infty} \le \frac{2}{|x|}$. \\
We now study the term $I_2$, which is $I_{2,1} + I_{2,2} := $
$$\int_{\left \{z : |z| > \frac{ |x|}{2 \left \| \gamma \right \|_\infty }, \\ |x + s z \gamma(x)| \le 1 \right \} } \int_0^1 (\gamma (x) \cdot z)^T \cdot H^2 f_{(x + s \gamma (x) \cdot z)} \cdot (\gamma (x) \cdot z) ds dz +$$
$$ + \int_{ \left \{z : |z| > \frac{ |x|}{2 \left \| \gamma \right \|_\infty }, \\ |x + s z \gamma(x)| > 1 \right \} } \int_0^1 (\gamma (x) \cdot z)^T \cdot H^2 f_{(x + s \gamma (x) \cdot z)} \cdot (\gamma (x) \cdot z) ds dz.$$
On $I_{2,1}$ we can upper bound the hessian matrix with $ c \epsilon$ and so we get 
\begin{equation}
I_{2,1} \le c \epsilon  \int_{\left \{z : |z| > \frac{ |x|}{2 \left \| \gamma \right \|_\infty }, \\ |x + s z \gamma(x)| \le 1 \right \} } |z|^2 \left \| \gamma \right \|_\infty^2 F(z) dz \le c \epsilon  \int_{\left \{z : |z| > \frac{ |x|}{2 \left \| \gamma \right \|_\infty } \right \} } |z|^2 e^{\epsilon  |z|} e^{ - \epsilon  |z|} F(z) dz \le  c \epsilon e^{ -\epsilon \frac{ |x|}{2 \left \| \gamma \right \|_\infty }} ,
\label{eq: I21}
\end{equation}
where in the last inequality we have used that $|z| > \frac{ |x|}{2 \left \| \gamma \right \|_\infty }$; after that we have enlarged the domain of integration and used the fifth point of Assumption A3 to upper bound the integral with a constant. \\
On $I_{2,2}$ we still use \eqref{eq: deriv seconde f*}, getting
$$I_{2,2} \le c \int_{ \left \{z : |z| > \frac{ |x|}{2 \left \| \gamma \right \|_\infty }, \\ |x + s z \gamma(x)| > 1 \right \} } \int_0^1 |z|^2 \left \| \gamma \right \|_\infty^2  \epsilon e^{\epsilon |x + s z \gamma(x)|} (\epsilon + \frac{1}{|x + s z \gamma(x)|}) F(z) dz ds  \le $$
\begin{equation}
\le c(\epsilon + 1) \epsilon e^{\epsilon |x|} \int_{ \left \{z : |z| > \frac{ |x|}{2 \left \| \gamma \right \|_\infty }, \\ |x + s z \gamma(x)| > 1 \right \} } |z|^2 e^{\epsilon |z| \left \| \gamma \right \|_\infty } F(z) dz.
\label{eq: I22}
\end{equation}
We define $J(x) := \int_{ \left \{z : |z| > \frac{ |x|}{2 \left \| \gamma \right \|_\infty } \right \} } |z|^2 e^{\epsilon |z| \left \| \gamma \right \|_\infty } F(z) dz$. From \eqref{eq: I1 f*} \eqref{eq: I21} and \eqref{eq: I22}, observing that the domain of integration of the integral defined as $J$ contains the one in \eqref{eq: I22} and using the boundedness of $\gamma$, we get
\begin{equation}
|A_d f^* (x)| \le c \, \epsilon e^{\epsilon |x| }(\epsilon + \frac{1}{|x|} + e^{- c \epsilon |x|} + (\epsilon + 1) J(x)). 
\label{eq: Ad f*}
\end{equation}
From \eqref{eq: Ac f*} and \eqref{eq: Ad f*} we get, using also the boundedness of $a$ and J which follows from A1 the fifth point of A3,
\begin{equation}
|A f^* (x)| \le \epsilon e^{\epsilon |x| } (c \epsilon - \tilde{C}) + c \epsilon e^{\epsilon |x| } (\frac{1}{|x|} + e^{- c \epsilon |x|} + J(x)).
\label{eq: final Af*}
\end{equation}
Since $\epsilon > 0$ can be chosen small; $c \epsilon $ is therefore less then $\tilde{C}$ and so the first term here above turns out being negative. \\
Moreover, the second term on the right hand side of \eqref{eq: final Af*} is $o(f^*)$. Indeed, $c \epsilon (\frac{1}{|x|} + e^{- c \epsilon |x|})$ clearly goes to $0$ for $|x| \rightarrow \infty$ and the domain of integration of the integral defined as $J(x)$  is the set of $z$ such that $|z| > \frac{ |x|}{2 \left \| \gamma \right \|_\infty } $. Therefore, for $|x| \rightarrow \infty$ the contribution of the integral becomes null. \\
It follows $A f^* \le - c_1 f^* + o (f^*)$, as we wanted. \\
We get the drift condition holding true on the proposed function $f^*$ and so the Assumption 3* does. 
\end{proof}

\subsection{Proof of Lemma \ref{lemma: convolution}}
\begin{proof}
We observe that, from the definition of convolution and Cauchy-Schwartz inequality it is
$$|f * g (x)| = \int_{\mathbb{R}^d} |f(x - y)| |g(y)| dy = \int_{\mathbb{R}^d} |f(x - y)|^{\frac{1}{2}} |g(y)| |f(x - y)|^{\frac{1}{2}} dy \le$$
$$ \le c (\int_{\mathbb{R}^d} |f(x - y)| |g(y)|^2 dy)^\frac{1}{2} (\int_{\mathbb{R}^d} |f(x - y)| dy)^\frac{1}{2} = c (\int_{\mathbb{R}^d} |f(x - y)| |g(y)|^2 dy)^\frac{1}{2} \left \|f \right \|_{1, \mathbb{R}^d}^\frac{1}{2}.$$
Therefore the following bound on the $L^2$ norm on $A$ holds true:
\begin{equation}
\left \| f*g \right \|_A^2 = \int_A |f*g (x)|^2 dx \le c \left \|f \right \|_{1, \mathbb{R}^d} \int_A (\int_{\mathbb{R}^d} |f(x - y)| |g(y)|^2 dy) dx .
\label{eq: convolution stima}
\end{equation}
We now use the fact that $\mathcal{S}$, the support of $f$, satisfies $diam(\mathcal{S}) < 2 \sqrt{d}$ and so, since $x \in A$, it follows $ \int_{\tilde{A}^c} |f(x - y)| |g(y)|^2 dy = 0$ for $\tilde{A}$ compact set of $\mathbb{R}^d$ such that $\tilde{A} := \left \{ \zeta \in \mathbb{R}^d : d(\zeta, A) \le 2 \sqrt{d} \right \}$. Using moreover Fubini's theorem, the right hand side of \eqref{eq: convolution stima} becomes
$$c \left \|f \right \|_{1, \mathbb{R}^d} \int_A (\int_{\tilde{A}} |f(x - y)| |g(y)|^2 dy) dx = c \left \|f \right \|_{1, \mathbb{R}^d} \int_{\tilde{A}} |g(y)|^2(\int_A |f(x - y)| dx) dy \le c \left \|f \right \|_{1, \mathbb{R}^d}^2 \left \|g \right \|_{2, \tilde{A}}^2, $$
where in the last estimation we have enlarged the integration domain of $f$ to $\mathbb{R}^d$.
The estimation here above joint with \eqref{eq: convolution stima} gives us 
$$\left \| f*g \right \|_A^2\le c \left \|f \right \|_{1, \mathbb{R}^d}^2 \left \|g \right \|_{2, \tilde{A}}^2.$$

\end{proof}


\begin{thebibliography}{4}


\bibitem{Applebaum} Applebaum, David. L\'evy processes and stochastic calculus. Cambridge university press, 2009.

\bibitem{BarShe} Barndorff-Nielsen, O. E. and Shephard, N. (2001). Non-Gaussian Ornstein-Uhlenbeck-based models and some of their uses in financial economics. J. R. Stat. Soc., Ser. B, Stat. Methodol., 63, 167-241.

\bibitem{Bates} Bates, D.S. (1996). Jumps and Stochastic Volatility: Exchange Rate Processes Implicit in Deutsche Mark. The Review of Financial Studies, 9(1), 69-107.

\bibitem{Chen} Chen, Z. Q., Hu, E., Xie, L., Zhang, X. (2017). Heat kernels for non-symmetric diffusion operators with jumps. Journal of Differential Equations, 263(10), 6576-6634.

\bibitem{Decomp} Comte, F., Lacour, C. (2013). Anisotropic adaptive kernel deconvolution. In Annales de l'IHP Probabilités et statistiques (Vol. 49, No. 2, pp. 569-609).

\bibitem{Main adapt} Comte, F., Prieur, C., Samson, A. (2017). Adaptive estimation for stochastic damping Hamiltonian systems under partial observation. Stochastic processes and their applications, 127(11), 3689-3718.

\bibitem{RD} Dalalyan, A. and Reiss, M. (2007). Asymptotic statistical equivalence for ergodic
diffusions: the multidimensional case. Probab. Theory Relat. Fields, 137(1), 25–47.

\bibitem{Neuro} Ditlevsen, S. and Greenwood, P. (2013). The Morris–Lecar neuron model embeds a leaky integrate-and-fire model. Journal of Mathematical Biology 67 239-259.

\bibitem{Mixing} Doukhan, P. (2012). Mixing: properties and examples (Vol. 85). Springer Science and Business Media.

\bibitem{GLM} Gloter, A., Loukianova, D. and Mai, H. (2018). Jump filtering and efficient drift estimation for L\'evy-driven SDEs. The Annals of Statistics, 46(4), 1445-1480.

\bibitem{Adaptive} Goldenshluger, A., Lepski, O. (2011). Bandwidth selection in kernel density estimation: oracle inequalities and adaptive minimax optimality. The Annals of Statistics, 39(3), 1608-1632.

\bibitem{KR} Klein, T., Rio, E. (2005). Concentration around the mean for maxima of empirical processes. The Annals of Probability, 33(3), 1060-1077.

\bibitem{Kou} Kou, S.G. (2002). A Jump-Diffusion Model for Option Pricing. Management Science, 48, 1086-1101.

\bibitem{Kus_Yos} Kusuoka, S., Yoshida, N. (2000). Malliavin calculus, geometric mixing, and expansion of diffusion functionals. Probability Theory and Related Fields, 116(4), 457-484.

\bibitem{Kut} Kutoyants, Y. A. (2013). Statistical inference for ergodic diffusion processes. Springer Science and Business Media.

\bibitem{Lacour et al} Lacour, C., Massart, P., Rivoirard, V. (2017). Estimator selection: a new method with applications to kernel density estimation. Sankhya A, 79(2), 298-335.

\bibitem{Lep} Lepski, O. (2013). Multivariate density estimation under sup-norm loss: oracle approach, adaptation and independence structure. The Annals of Statistics, 41(2), 1005-1034.

\bibitem{Lep99} Lepski, O. V., Levit, B. Y. (1999). Adaptive non-parametric estimation of smooth multivariate functions.

\bibitem{18 GLM} Masuda, H. (2007). Ergodicity and exponential beta - mixing bounds for multidimensional diffusions with jumps. Stochastic processes and their applications, 117(1), 35-56.

\bibitem{Merton} Merton, R.C. (1976). Option pricing when underlying stock returns are discontinuous. Journal of Financial Economics, 3, 125-144

\bibitem{May_Twe} Meyn, S. P., Tweedie, R. L. (1993). Stability of Markovian processes III: Foster–Lyapunov criteria for continuous-time processes. Advances in Applied Probability, 25(3), 518-548.

\bibitem{Nik} Nikolskii, S. M. Approximation of Functions of Several Variables and Embedding Theorems (Springer, Berlin, 1975).

\bibitem{Str_Twe} Stramer, O., Tweedie, R. L. (1997). Existence and stability of weak solutions to stochastic differential equations with non-smooth coefficients. Statistica Sinica, 577-593.

\bibitem{Strauch} Strauch, C. (2018). Adaptive invariant density estimation for ergodic diffusions over anisotropic classes. The Annals of Statistics, 46(6B), 3451-3480.

\bibitem{Ver} Veretennikov, A. Y. (1988). Bounds for the mixing rate in the theory of stochastic equations. Theory of Probability and Its Applications, 32(2), 273-281.
ISO 690	

\bibitem{Viennet} Viennet, G. (1997). Inequalities for absolutely regular sequences: application to density estimation. Probability theory and related fields, 107(4), 467-492.

\end{thebibliography}
\end{document}